\documentclass[12pt,twoside,english]{article}
\usepackage{amssymb,amsmath,babel,}

\newtheorem{Lemma}{Lemma}[section]
\newtheorem{Theorem}[Lemma]{Theorem}

\newtheorem{Proposition}[Lemma]{Proposition}

\newtheorem{Definition}[Lemma]{Definition}
\newtheorem{Remark}[Lemma]{Remark}
\newtheorem{Conjecture}[Lemma]{Conjecture}

\newcommand{\gp}{\mathfrak{p}}
\newcommand{\ZZ}{\mathbb{Z}}

\newcommand{\FF}{\mathbb{F}}
\newcommand{\CC}{\mathbb{C}}
\newcommand{\QQ}{\mathbb{Q}}

\newcommand{\NN}{\mathbb{N}}

\newcommand{\partialbis}{\partial^{(E)}}

\newcommand{\ie}{\emph{i.e.}}

\newcommand{\eg}{\emph{e.g.}}

\newcommand{\sqm}[4]{
\left(\begin{array}{ll}#1 & #2 \\ #3 & #4\end{array}\right)}
\newcommand{\binomial}[2] {{\binom{#1}{#2}}} 
\newcommand\CVD{{\hfill\hfil{\lower 2 pt\hbox{\vrule\vbox to 7pt 
{\hrule width 6pt\vfill\hrule}\vrule}}}\vskip 0.5cm}

\let\le=\leq  
\let\ge=\geq

\title{On certain families of \\ 
Drinfeld quasi-modular forms}
\author{Vincent Bosser, Federico Pellarin}

\begin{document}

\maketitle

\tableofcontents

\section{Introduction.}


The primary aim of this paper is to describe some partial advances in the solution of the following two problems.
\begin{enumerate}
\item Find the maximal {\em order of vanishing at infinity} of a non-zero {\em Drinfeld quasi-modular form}
of given weight.
\item Determine {\em differential properties} of Drinfeld quasi-modular forms of given weight and depth
with maximal order of vanishing at infinity (these forms will be called {\em extremal}). 
\end{enumerate}

Our results are obtained in a constructive way, studying families of forms with peculiar properties.
For our purposes, general tools need to be developed. Some will appear of independent interest.

Before going deeper in these topics and rigorously defining the entities
above, we present an overview of the more
familiar framework of quasi-modular forms on the complex upper-half
plane, for the group $\mathbf{SL}_2(\ZZ)$.

\subsection{The classical framework.\label{section:classical}}
Let $z=x+\boldsymbol{i}y\in\CC$ with $y>0$ and $u\in\CC$ be complex
numbers related by $u=e^{2\pi \boldsymbol{i} z}$, so that $0<|u|<1$.
For $i=1,2,3$ the series:
\begin{equation*}
E_{2i}(z) =  1+b_i\sum_{n=1}^\infty n^{2i-1}\frac{u^n}{1-u^n},
\end{equation*}
with $b_1=-24,b_2=240,b_3=-504$, normally converge in any compact subset of the domain determined by $|u|<1$
and represent algebraically independent functions.
The $\CC$-algebra of {\em quasi-modular forms} is the polynomial algebra $\widetilde{M}:=\CC[E_2,E_4,E_6]$,
which is graded by the {\em weights} (where the weight of $E_{2i}$ is
$2i$ for $i=1,2,3$) and filtered by the {\em depths} (the depth of a
polynomial in $\widetilde{M}$ is its degree in $E_2$):
$$\widetilde{M}=\bigoplus_{w\geq 0}\bigcup_{l\geq 0}\widetilde{M}_w^{\leq l},$$
where $\widetilde{M}_w^{\leq l}$ is the $\CC$-vector space spanned by the forms
of weight $w$ and {\em depth} $\leq l$.

Any element $f$ of $\widetilde{M}_w^{\leq l}\setminus\{0\}$ has a non-vanishing $u$-expansion
\begin{equation}\label{expansion}
f(z)=\sum_{i=0}^{\infty}c_iu^i,\quad c_i\in\CC
\end{equation}
and the natural problem of determining the image of the function \begin{equation}\label{nuinfty}
\begin{array}{rcl} \nu_\infty:\widetilde{M}^{\leq l}_w\setminus\{0\} & \rightarrow & \ZZ_{\geq 0} \\
f & \mapsto & \inf\{i,c_i\neq 0\},\text{ with $f$ as in
(\ref{expansion})}\end{array}\end{equation} arises,
for given $l,w$. 

In \cite[p. 459]{KK1}, Kaneko and Koike 
ask whether the image of $\nu_\infty$ on the set $\widetilde{M}^{\leq l}_w\setminus\{0\}$
is precisely the interval $[0,1,\ldots,\dim_{\CC}(\widetilde{M}^{\leq l}_w)-1]$, as numerical
investigations suggest, for $w$ small. This property, if true, would  imply that for any $f\in\widetilde{M}^{\leq l}_w\setminus\{0\}$,
\begin{equation}\label{eq:estimateKK}\nu_\infty(f)\leq
\dim_{\CC}(\widetilde{M}^{\leq l}_w)-1
\leq \frac{1}{12}(w+11l+l(w-l)).\end{equation}

We recall that $\widetilde{M}$ is a $D$-differential algebra, with
$D:=u\frac{d}{du}=(2\pi\boldsymbol{i})^{-1}\frac{d}{dz}$.
Looking at the resultant $\mathbf{Res}_{E_2}(f,Df)$ of $f$ and $Df$,
seen as polynomials in $E_2$, it is not difficult
to prove that for $f$ as above, irreducible (\footnote{The inequality
(\ref{eq:saradha}) also holds for $f$ not necessarily irreducible.}),
\begin{equation}\label{eq:saradha}
\nu_\infty(f)\leq \frac{1}{12}(w+2l(w-l)).
\end{equation}
Similar inequalities have already been used to describe diophantine properties of
certain complex numbers, see \eg\ \cite{grinspan, saradha}.
%
%
%
%
In fact, in order to prove (\ref{eq:saradha}), the above resultant
can only be used if it does not vanish, that is, if $Df$ is not divisible
by $f$. To prove the estimate for the remaining forms,
we need to characterise those forms $f$ such that $f$ divides $Df$.
The key point is here to remark that the only principal
prime ideal of $\widetilde{M}$ which is stable
by differentiation is the ideal $(\Delta)$, where $\Delta$
is the discriminant function (see \cite[Chapter 10, Lemma 5.2]{NP}).
The remaining case $f=\Delta$ in (\ref{eq:saradha}) can then be
checked directly.

Apart from some choices $(l,w)$ with $1\leq l\leq 10$ and $w\leq 20$
and the case $l=0$, the upper bound of (\ref{eq:saradha})
is weaker than that of (\ref{eq:estimateKK}).
The truth of the sharper estimate (\ref{eq:estimateKK}) remains
unknown for general $l,w$.

Let  $l,w$ be integers such that $\widetilde{M}^{\leq l}_w\not=(0)$, and let 
$f_{l,w}$  be the unique non-vanishing normalised
(\footnote{A formal series $f=\sum_{i\geq i_0}c_iu^i$
with $c_{i_0}\neq 0$ is said to be {\em normalised}
if $c_{i_0}=1$. A quasi-modular form is normalised if its
$u$-expansion is.}) quasi-modular form
of the space $\widetilde{M}^{\leq l}_w$
with the property that the function (\ref{nuinfty}) attains its
maximal value on it.

In \cite[Theorem 2.1]{KK1}, Kaneko and Koike constructed a family of
quasi-modular forms which turns out to be, up to a non-zero scalar
factor and by means of elementary arguments, the family
$(f_{1,2i})$ (\footnote{In \cite{KK1}, Kaneko and Koike call any
non-vanishing form
$f\in\widetilde{M}^{\leq l}_k\setminus\widetilde{M}^{\leq l-1}_k$
for which the equality
$\nu_{\infty}(f)=\dim_{\CC}(\widetilde{M}^{\leq l}_k)-1$ holds,
{\em extremal}. We warn the reader
that in this paper, we will use this terminology in a different way.}).
They prove that
$\nu_\infty(f_{1,2i})=
\lfloor\frac{i}{3}\rfloor=\dim(\widetilde{M}^{\leq 1}_{2i})-1$ ($i\geq 0$),
where $\lfloor\cdot\rfloor$ denotes the lower integer part of a real number.
Their construction
is performed with an inductive process, in which
the differential operators (``Serre's differential operators",
cf. p. 467 of \cite{KK1}) defined by:
\begin{eqnarray}
\theta_d^{(n-1)}f&=&D^nf-\frac{n+d-1}{12}[E_2,f]^{(2,d)}_{n-1},\label{eq:rcb}
\end{eqnarray}
play a crucial role,
$[f,g]^{(2,d)}_{n-1}$ being a suitable normalisation of the {\em Rankin-Cohen bracket} of $f$ and
$g$, of order $n-1$ and weights $2$ and $d$ (its definition is available, for example, on p. 466 of
\cite{KK1}). It turns out that for all $i$, the form $f_{1,2i}$ is essentially unique satisfying: 
\begin{equation}
\theta_{2i-1}^{(1)}f_{1,2i}=\lambda_{1,2i}\Delta f_{1,2i-8},\quad i\geq 1,\label{eq:fstar}
\end{equation}
where $\Delta=(E_4^3-E_6^2)/1728$ and $\lambda_{1,2i}\in\QQ$ with $\lambda_{1,2i}=0$ if and only if $3|i$.

%

For $l=2$, Kaneko and Koike develop similar constructions in \cite{KK1}. We omit to describe their results 
referring to \cite{KK2, KK1} for further details and references; however, we did not find any other result in the direction of inequalities (\ref{eq:estimateKK})
for $l\geq 2$.

\subsection{The drinfeldian framework: our results.}

From now on, the symbols $u,\widetilde{M},\Delta,\nu_\infty,D$
will be used with a new meaning which we now describe. Occasionally, the older meanings
related to the classical framework will reappear, but the reader should
not encounter any trouble with these occurrences.

In the following, $q=p^e$ is a positive power of a given prime number
$p$ fixed once and for all, and $\theta$
will be an indeterminate over $\FF_q$. Certain results of this text
do not hold for certain choices of $q$; this will be highlighted on a case to case basis. 

%

Let $C$ be the completion of an algebraic closure $\overline{K_\infty}$
of the field $K_\infty:=\FF_q((1/\theta))$
for the unique extension of the valuation $-\deg_\theta$ to $\overline{K_\infty}$ (the $\theta^{-1}$-adic
valuation). For this valuation, $K_\infty$ is the completion of 
its subfield $K=\FF_q(\theta)$. We will also denote by $A$ the $\FF_q$-algebra $\FF_q[\theta]$.
On $C,K_\infty,K,A$, we will denote by $|\cdot|$ the ultrametric absolute value $q^{\deg_\theta(\cdot)}$.

The story of {\em Drinfeld modular forms} (\footnote{Rigid analytic
modular forms, where the ``base fields" are global, of positive
characteristic.})
begins with the pioneering work of Goss \cite{Go1,Go2}. Later, in the
important paper \cite{Ge}, Gekeler considered the two functions
$E,g$ and discovered the function $h$, allowing us, later in
\cite{BP}, to investigate some properties of the $C$-algebra 
of {\em Drinfeld quasi-modular forms} $\widetilde{M}:=C[E,g,h]$
%
(\footnote{Notice that in the introduction of \cite{BP}, the sentence
attributing to Gekeler's paper the first occurrence of modular forms for
$\mathbf{GL}_2(\FF_q[\theta])$ is obviously incorrect.}). As ``drinfeldian framework" (title of this section)
we mean a natural counterpart of the theory sketched in Section \ref{section:classical}, for these
Drinfeld quasi-modular forms.

We refer to the first part of our paper \cite{BP} for all the basic properties of Drinfeld quasi-modular forms, noticing that the indeterminate
$T$ there becomes $\theta$ here (\footnote{Several notations of \cite{BP} change in the present work. We  
adopt at the same time notations issued from Gekeler as in \cite{Ge} and notations that will be compatible with 
Papanikolas \cite{Pa} since we believe that the use of $t$-motives will eventually intervene in the
theory of Drinfeld modular forms, and we want to keep certain symbols, such as $t$, free for that occasion.}).
Gekeler's functions $E,g,h$ are algebraically independent quasi-modular forms for the
(homographical) action of the group $\mathbf{GL}_2(A)$ on $\Omega:=C\setminus K_\infty$.
For the three functions $E$, $g$, $h$,
the triples $(w,m,l)\in\ZZ_{\geq 0}\times(\ZZ/(q-1)\ZZ)\times
\ZZ_{\geq 0}$,  where $w$ is the
weight, $m$ the type and $l$ the depth, are $(2,1,1)$,
$(q-1,0,0)$ and $(q+1,1,0)$, respectively. 

If we denote by $\widetilde{M}^{\leq l}_{w,m}$ the $C$-vector space of Drinfeld quasi-modular forms of weight $w$, type $m$ and depth $\leq l$
(which is by definition the space $(0)$ if $l<0$),
we have \cite[Proposition 2.2]{BP}:
$$ \widetilde{M}=
\bigoplus_{\genfrac{}{}{0pt}{1}{w\in\ZZ_{\geq 0}}{m\in\ZZ/(q-1)\ZZ}}\bigcup_{l\geq 0}
\widetilde{M}^{\leq l}_{w,m}.$$

\subsubsection{First part: extremality}

We recall \cite[Section 2 and Lemma 4.2 (i)-(iii)]{BP} that $E,g,h$ have $u$-expansions 
convergent in a neighbourhood of $u=0$,
with $u(z)=1/e_{\mathcal{C}}(\widetilde{\pi}z)$, where $e_{\mathcal{C}}$ is Carlitz's exponential
function and $\widetilde{\pi}\in C$
is one of its fundamental periods (chosen once and for all) (\footnote{In \cite{BP}, we wrote about
$t$-expansions instead of $u$-expansions, and the period of Carlitz's exponential was denoted by
$\overline{\pi}$ instead of $\widetilde{\pi}$, according to Gekeler \cite{Ge}.}). There is a $C$-algebra
homomorphism $\widetilde{M}\hookrightarrow C[[u]]$. We will identify $f\in\widetilde{M}$
with its image $\sum_{i\geq 0}c_iu^i$ in $C[[u]]$.  We write $\nu_\infty(0)=\infty$ and, for
$f\neq 0$, $\nu_\infty(f)= \min\{i\text{ such that }c_i\neq 0\}$.

\begin{Definition}\emph{
Let $w,l$ be fixed non-negative integers and $m$ be a class in $\ZZ/(q-1)\ZZ$.
A non-zero quasi-modular form $f\in\widetilde{M}^{\leq l}_{w,m}$ is
said to be {\em extremal of depth $\leq l$} if for all
$g\in\widetilde{M}^{\leq l}_{w,m}\setminus\{0\}$, $\nu_\infty(g)\leq \nu_\infty(f)$.}
\end{Definition}


A quasi-modular form which is extremal of depth $\leq l$ needs not to
be extremal of depth $\leq l'$ for $l'>l$.
On the other hand, a straightforward argument shows that, for $w,m, l$
fixed, a quasi-modular form of weight $w$ and type $m$
which is extremal of depth $\leq l$ is uniquely determined, up to
multiplication by an element of $C^\times$.

When a quasi-modular form $f\in\widetilde{M}^{\leq l}_{w,m}\setminus\widetilde{M}^{\leq l-1}_{w,m}$
is extremal of depth $\leq l$, we will often say that $f$ is {\em extremal}. We will adopt this simplified terminology when
the context allows the complete determination of $l$.

For all $l,w,m$ such that $\widetilde{M}^{\leq l}_{w,m}\neq(0)$, let
$f_{l,w,m}\in\widetilde{M}^{\leq l}_{w,m}$ be the unique normalised extremal quasi-modular
form of depth $\leq l$.
Contrary to the classical framework, $u$-expansions of Drinfeld quasi-modular forms
are difficult to compute, and Gekeler's algorithms developed in \cite{Ge} are required. Thanks to them,
we did experiments that, after observation of the cases $q=2,3,5$, $w\leq q^3+1$,
$l\leq q^2+1$ and any value of $m$, suggest the existence of a (conjectural) estimate as follows. For all $\epsilon>0$ and for all
$l$ big enough depending on $\epsilon$:
\begin{equation}
\nu_\infty(f_{l,w,m})\leq (1+\epsilon)l(w-l)
\end{equation}
(notice that if $l>0$, $w>l$).

Just like inequality (\ref{eq:estimateKK}), this inequality seems to be
rather difficult to prove. Even weaker estimates like an analogue of
(\ref{eq:saradha}) are presently unavailable (see discussion in
Section 1 of \cite{BP}). This is essentially because
in our case there are infinitely many irreducible quasi-modular forms
$f$ such that $Df$ is divisible by $f$, and in this last case there
is no obvious candidate to replace the modular form $\mathbf{Res}_{E}(f,Df)$. 

It would be interesting to find an explicit function $c(w,l)$ of the weight and the depth such that for all $l\gg 0$,
$\nu_\infty(f_{l,w,m})\leq c(w,l)$. Showing 
the existence of a constant $c(q)>0$, depending on $q$ only, such that if $l>0$,
$$\nu_\infty(f_{l,w,m})\leq c(q)l(w-l),$$ would also have interesting arithmetic consequences.
For instance, the results of the present paper 
show that if $c(q)$ exists, then $c(q)\geq 1$ (cf. Proposition \ref{x0x1x2}).
In the first part of this paper,
we discuss partial advances toward these estimates for small depths.

Although analogues of ``higher Serre's operators'' can be constructed
(this paper, Section \ref{sserre_operators}), ideas of proof of Kaneko and Koike as in \cite{KK2} cannot extend to our case because these operators have too big kernels due to the positive characteristic
(but see Section \ref{hypergeometric}
for some condition analogous to (\ref{eq:fstar})).

In section \ref{extremality} we study the sequence of Drinfeld quasi-modular forms $(x_k)_{k\geq 0}$, with 
$x_k\in\widetilde{M}_{q^k+1,1}^{\leq 1}\setminus \widetilde{M}_{q^k+1,1}^{\leq 0}$, defined by $x_0=-E$, $x_1=-Eg-h$ and by the recursion formula
\[x_k=x_{k-1}g^{q^{k-1}}-[k-1]x_{k-2}\Delta^{q^{k-2}},\quad k\geq 2,\] where $\Delta:=-h^{q-1}$.
After having dealt with some basic properties of this sequence, 
we show in Proposition \ref{xkordqk} that, for all $k\geq 0$, $\nu_\infty(x_k)=q^k$.
This shows that in general,
$$\nu_{\infty}(f_{l,w,m})>\dim_C(\widetilde{M}^{\leq l}_{w,m})-1,$$ in apparent contradiction with Kaneko and Koike's 
prediction (\ref{eq:estimateKK}) in the classical framework. 

We can show that for all $k$, $x_k$ is an irreducible polynomial in $E,g,h$ and a resultant argument (see Section \ref{extremality_xy})
yields:
\begin{Theorem}\label{theo_xk_extremal} For all $k\geq 0$, $x_k$ is an extremal quasi-modular form.\end{Theorem}
Looking back at subsequences of Kaneko and Koike's sequence $(f_{1,2i})_{i\geq 0}$, there does not seem to be similar recursion
formulas, with weights varying as sequences like $(\alpha q^k+\beta)_{k\geq 0}$ rather than as arithmetic progressions.
But experimentally, congruences between $u$-expansions of certain forms $f_{1,2i}$'s seem to occur. They could 
be consequence of Clausen-von Staudt congruences for Bernoulli numbers.

Our investigation was pushed a step further, with the sequence $(\xi_k)_{k\geq 0}$ defined, for $k\geq 0$ by:
\begin{equation}\label{def_xi_k}
\xi_k=[k]^qx_{k+1}x_{k-1}^q-[k+1]x_k^{q+1}\in\widetilde{M}_{(q+1)(q^k+1),2}^{\leq q+1}\setminus\widetilde{M}_{(q+1)(q^k+1),2}^{\leq q},\end{equation}
where $[k]:=\theta^{q^k}-\theta$ ($k\geq 1$), $[0]:=1$, and we have set $x_{-1}:=-h^{1/q}$. Again, we could compute $\nu_\infty(\xi_k)$ for all $k$ and prove that $\xi_k$
is always irreducible, implying the following (Section \ref{sectionxik}):
\begin{Theorem}\label{xik_extremal} Assuming that $q\geq 3$,
the form $\xi_k$ is extremal for all $k\geq 0$.
\end{Theorem}
As a product of these investigations, we obtain the following multiplicity
estimate:
\begin{Theorem}\label{th_multiplicity3} Let $w$ and $m$ be integers
such that $\widetilde{M}^{\leq q^2}_{w,m}\not=\{0\}$, and let $f$ be a
non-vanishing form in $\widetilde{M}^{\leq q^2}_{w,m}$. Then
$$\nu_\infty(f)\leq (q^3+1)(w-l).$$
\end{Theorem}

\subsubsection{Second part: differential properties.}

We were surprised to remark that the forms $x_k$ and $\xi_k$ also enjoy a rich differential structure, and the
second part of this text, Section \ref{differential_extremality}, is devoted to reporting our knowledge
on this topic. 

In all the following, we write $D=(D_n)_{n\geq 0}$ for the collection of higher derivatives on the
$C$-algebra of holomorphic functions on $\Omega=C\setminus K_\infty$ introduced in \cite[Section 1]{BP}.
Therefore, $D_1=(-\widetilde{\pi})^{-1}d/dz=u^2d/du$.
By Theorem 2 of loc. cit., $D_n$ induces a $C$-linear map $$\widetilde{M}_{w,m}^{\leq
l}\xrightarrow{D_n}\widetilde{M}_{w+2n,m+n}^{\leq l+n},$$
so that the $C$-algebra $\widetilde{M}$ is {\em $D$-stable} (or {\em hyperdifferential}).

Already in \cite{BP}, we have remarked that the problem of estimating the quantity $\nu_\infty(f)$
for a Drinfeld quasi-modular form $f$ is intimately related to its differential properties
(this point of view was inherited by Nesterenko,
and finds its foundations in Siegel and Shidlowski's work). 

In the papers \cite{Go1,Go2}, $u$-expansions were already considered, and their behaviour 
immediately appeared to be surprisingly erratic.
Later, in \cite{Ge}, Gekeler described algorithms to compute their $u$-expansions.
However, the unpredictable character of the
coefficients of the $u$-expansions of all these Drinfeldian forms remains nowadays one of the typical aspects making this theory 
independent from the classical one. Similar observations can be made concerning the problem of Hecke's eigenforms.
As a strange resurgence of older problems, we remarked in \cite{BP} that also the operators $D_n$ behave erratically over the polynomial algebra
$C[E,g,h]$.

Here, we introduce the following:


\begin{Definition}{\em
Let $f$ be a non-zero element of $\widetilde{M}$. We define the {\em differential exponent}
$\epsilon_D(f)$ of $f$ as follows: it is $\infty$ if $\frac{D_{n}f}{f}\in\widetilde{M}$ for all $n\geq 1$,
and otherwise it is the smallest integer $k\geq 0$ satisfying $\frac{D_{p^k}f}{f}\not\in\widetilde{M}$
(thus $\epsilon(f)=0$ means $\frac{D_{1}f}{f}\not\in\widetilde{M}$).

Let $l,w,m$ be such that
$\widetilde{M}^{\leq l}_{w,m}\setminus\widetilde{M}^{\leq l-1}_{w,m}\neq\emptyset$
(with the convention that if $l<0$, $\widetilde{M}_{w,m}^{\leq l}:=(0)$), let 
$f$ be a quasi-modular form in this set. We say that $f$ is {\em differentially 
extremal} of weight $w$, type $m$, depth $l$, if it attains the biggest {\em finite} differential exponent within this set.
}\end{Definition}
By \cite[Proposition 3.6]{BP}, a differentially extremal quasi-modular form cannot be proportional to a power of $h$ and 
if $q\neq 2,3$ we obtained, in \cite[Theorem 3]{BP}, that if $f$ is not proportional to a power of $h$, then $f$ has a 
finite differential exponent.
In Section \ref{differential_extremality}, by using a result of Cornelissen in \cite{GG}
on the factorisation in $K[g,h]$ of certain normalised Eisenstein series, we prove:
\begin{Theorem}\label{theorem:xkbis} For all $k\geq 0$,
$x_k$ is differentially extremal of differential exponent
$(k+1)e$.\end{Theorem}
We are presently unable to show the differential extremality of the forms $\xi_k$, but in
Section \ref{experiments}, we describe numerical computations which seem to confirm this hypothesis.

Sections \ref{hypergeometric} and \ref{Fk} develop two questions, partly independent but not completely disjoint to
the problem of finding differentially extremal quasi-modular forms. The content of these sections will be 
presented at their respective beginning; the reader can skip them in a first reading of the paper.

\subsubsection{Third part: differential tools.}
The proofs of the statements above require several technical tools which appear in Section \ref{differential_tools}
of this paper.
In this section, the reader can find several results, some of which are of
independent interest, described in the summary below.

Let $n,d$ be non-negative integers. We define the {\em $n$-th Serre's
operator of degree $d$}, $\partial_n^{(d)}:\widetilde{M}\rightarrow
\widetilde{M}$, by the formula
\begin{equation}\label{defpartial}
\partial_n^{(d)}f = D_nf + \sum_{i=1}^n (-1)^i{d+n-1\choose i}(D_{n-i}f)(D_{i-1}E).
\end{equation}
These operators can be considered as analogues of higher Serre's $C$-linear differential
operators (\ref{eq:rcb}) in the drinfeldian framework. In Theorem \ref{proposition_partial} of
Section~\ref{sserre_operators} we show the (not obvious at all) property that
$$\partial_n^{(w-l)}:\widetilde{M}^{\leq l}_{w,m}\rightarrow\widetilde{M}^{\leq l}_{w+2n,m+n}$$
(compare with \cite[Proposition 3.3]{KK1}).
The properties of these operators are essential in the proof of Theorem \ref{theorem:xkbis}.
A further application of the operators $\partial^{(d)}_n$ is contained in Section \ref{heckeoperators},
where we indicate a new technique to determine modular eigenforms of all the Hecke operators.

In Section \ref{algorithm} we furnish algorithms to compute the polynomials $$D_nE,D_ng,D_nh\in C[E,g,h].$$
These algorithms can be viewed as variants of Gekeler's algorithms in \cite{Ge}. Proposition \ref{propformulasgsxs}
is crucial, for example, in the computations of the polynomials $D_n\xi_k$ we made, as well as in the proof of 
Theorem \ref{theorem:xkbis}.

\subsubsection{Final remarks}

It is strongly possible that the forms $\xi_k$ are all differentially extremal. This can be checked, in principle,
by using the tools of Section \ref{algorithm}, up to tremendous calculations  we had not the courage to do,
but that can be done.
In fact, we got interested in these forms $\xi_k$ in an attempt of finding differentially extremal 
forms by solving linear equations; later, we found that these forms are extremal, yielding the actual 
presentation of this paper. This convinced us to 
follow a constructive approach to produce multiplicity estimates. 

 The reader may remark that the problem of finding such families is essentially of a diophantine nature. It can be proved, just as in \cite[p. 153]{KK2}, that 
the sequence $(x_k)_{k\geq 0}$ is related to the convergents of the continued fraction
expansion of the function $h/E$ as a formal series in $j:=g^{q+1}/\Delta$. On another side, in \cite{KK3},
Kaneko reveals a connection between certain modular solutions of the differential equations $\theta^{(1)}_{k}f=0$
and Apery's approximations of $\zeta(2)=\pi^2/6$.
Hence, it could reveal difficulty to explicitly construct new interesting families 
in higher depth. Nevertheless, we think that the connection between extremality and differential extremality of certain
families of Drinfeld quasi-modular forms is such a topic that will deserve further surprises.

\section{Order of vanishing and extremality.\label{extremality}}

As already mentioned in the introduction, the
main objective of this section is to introduce and study two families of extremal Drinfeld quasi-modular
forms: one in depth $\leq 1$ (the forms $x_k$, see Section~\ref{sectionxk})
and one in depth $\leq q+1$ (the forms $\xi_k$, see Section~\ref{sectionxik}).
We use these forms to prove multiplicity estimates
for quasi-modular forms of depth $\le q$ or $\leq q^2$
(Sections~\ref{sectionmultxk} and \ref{sectionmultxik}). The tools developed in this section allow to compute,
in Proposition \ref{x0x1x2}, certain extremal forms of depth $<\frac{q+1}{2}$.

\subsection{The family $(x_k)_{k\geq 0}$.\label{sectionxk}}

We begin by defining three sequences of Drinfeld quasi-modular forms:
$$(g_k)_{k\geq 0},(h_k)_{k\geq 0},(x_k)_{k\geq 0}.$$

Einsenstein's series for the group $\mathbf{GL}_2(A)$ are defined on p. 681 of
\cite{Ge}:
\begin{equation}\label{eq:definition_eisenstein}
E^{(w)}(z)= \sum_{a,b\in A}{}'(az+b)^{-w},
\end{equation}
where the dash 
indicates that the sum is restricted to $a,b$ not simultaneously zero.
It is easy to prove that the series $E^{(w)}$ converges uniformly on every
compact subset of $\Omega$, for all integers $w>0$, to  a Drinfeld modular
form of weight  $w$ and type $0$ which is non-zero if and only if
$w\equiv 0\pmod{q-1}$ (see \cite{Ge} p. 682).

Following \cite{Ge} p. 684 (and the notations of this reference),
let us write $g_0=1$ and, for $k\geq 1$,
\begin{equation}
g_k=(-1)^{k+1}\widetilde{\pi}^{1-q^k}L_kE^{(q^k-1)},
\label{def_gk}
\end{equation}
where $L_k:=[k][k-1]\cdots[1]$.

For $k\geq 0$, $g_k$ is a non-vanishing normalised modular form of weight $q^k-1$
and type $0$, whose expansion at infinity belongs to $A[[u]]$.

We have \cite[Proposition 6.9]{Ge}: $g_0=1$, $g_1=g$, and
\begin{equation}\label{eq:g0g1gk}
g_k=g_{k-1}g^{q^{k-1}}-[k-1]g_{k-2}
\Delta^{q^{k-2}},\quad k\ge 2.
\end{equation}

In order to define the forms $h_k$'s ($k\ge 0$) we introduce,
for fixed $w\in\NN$ and $m\in\ZZ/(q-1)\ZZ$,
Serre's operator $\partial_1^{(w)}:M_{w,m}\rightarrow\widetilde{M}$, which
is defined by
\begin{equation}\label{defserreop}
\partial_1^{(w)}f=D_1f-wEf.
\end{equation}
It is well-known (see \cite[Section (8.5)]{Ge}) that
\begin{equation*}
\partial_1^{(w)}(M_{w,m})\subset M_{w+2,m+1},
\end{equation*}
so we have in fact an operator $\partial_1^{(w)}:M_{w,m}\rightarrow
M_{w+2,m+1}$ (note that in \cite{Ge} Serre's operator is denoted by
$\partial_w$ and defined by the formula
$\partial_wf=wEf-D_1f=-\partial_1^{(w)}f$). 

We now define, for $k\geq 0$:
$$h_k=-\partial_1^{(q^k-1)}g_k,\quad k\geq 0.$$
For all $k$, $h_k$ is a modular form of weight $q^k+1$, type $1$. Moreover, $h_0=0$ and $h_1=h$ \cite[Theorem (9.1)]{Ge}.
Finally, we define the sequence $(x_k)_{k\geq 0}$ by:
$$x_0=-E\quad\text{\rm and}\quad x_k=D_1g_k,\quad k\geq 1.$$
We will see in a little while that this definition is compatible with that of the introduction.

Since by definition we have, for $k\ge 1$, $\partial_1^{(q^k-1)}f=D_1f+Ef$,
we find
$h_k=-D_1g_k-Eg_k=-x_k-Eg_k$. Hence the following identity holds
(one immediately checks that it is also true for $k=0$):
\begin{equation}\label{eq:egkhk}x_k=-Eg_k-h_k,\quad k\geq 0.\end{equation}
Therefore, the form $x_k$ is, for $k\geq 0$, a non-modular quasi-modular form of weight $q^k+1$, type $1$ and 
depth $1$.

It turns out that the three families $(g_k)_k$, $(h_k)_k$ and $(x_k)_k$
satisfy the same recursion formula.

\begin{Proposition}\label{sequences}
The sequence $(h_k)_{k\geq 0}$ is determined by the initial conditions
$$h_0=0,\quad h_1=h$$
and the recursion formula
\[h_k=h_{k-1}g^{q^{k-1}}-[k-1]h_{k-2}\Delta^{q^{k-2}},\quad k\geq 2.\]
Similarly, the 
sequence $(x_k)_{k\geq 0}$ is determined by the initial conditions
$$x_0=-E,\quad x_1=-Eg-h$$
and the recursion formula
\[x_k=x_{k-1}g^{q^{k-1}}-[k-1]x_{k-2}\Delta^{q^{k-2}},\quad k\geq 2.\]
\end{Proposition}


\medskip

\noindent {\em Proof.}
We begin with the recursion relation for $x_k$.
The formulas $x_0=-E$ (definition) and $x_1=D_1g=-Eg_1-h_1=-Eg-h$
have been already remarked.
If $k=2$, then $g_2=[1]h^{q-1}+g^{q+1}$, so that by formulas (2) of \cite{BP}:
\begin{eqnarray*}
x_2=D_1g_2&=&-[1]h^{q-1}E-(Eg+h)g^q\\
&=&x_1g^q-[1]x_0\Delta,
\end{eqnarray*}
which is the expected relation.
If now $k>2$, then $D_1g^{q^{k-1}}=0$ and $D_1\Delta^{q^{k-2}}=0$, so that,
by (\ref{eq:g0g1gk}):
\begin{eqnarray*}
x_k=D_1g_k&=&D_1(g_{k-1}g^{q^{k-1}})+D_1(-[k-1]g_{k-2}\Delta^{q^{k-2}})\\
&=&(D_1g_{k-1})g^{q^{k-1}}-[k-1](D_1g_{k-2})\Delta^{q^{k-2}}\\
&=&x_{k-1}g^{q^{k-1}}-[k-1]x_{k-2}\Delta^{q^{k-2}}. 
\end{eqnarray*}
The proof of the statement about the sequence $(h_k)_{k\geq 0}$ is now clear.
Indeed, by (\ref{eq:egkhk}) and the result on the sequence $(x_k)_{k\geq 0}$
we have just proved, we have, for $k\geq 2$:
$$-(Eg_k+h_k)=(-Eg_{k-1}-h_{k-1})g^{q^{k-1}}+[k-1](Eg_{k-2}+h_{k-2})\Delta^{q^{k-2}}.$$
Using now the recursion formula (\ref{eq:g0g1gk}) for the sequence
$(g_k)_{k\ge 0}$, we get the same recursion formula for $(h_k)_{k\ge 0}$.
\CVD

\subsubsection{Order of vanishing of the form $x_k$.}\label{sectionvanishing}

In this section, we determine the order of vanishing at infinity of the
form $x_k$ (for all $k$).

First of all, we recall our conventions for binomial cofficients
\cite[Section 3]{BP}. For
$n\in\ZZ$ and $i\in\ZZ$ with $i\geq 0$:
$${n\choose i}:=\frac{\prod_{k=1}^i (n-k+1)}{i!}. $$

We begin with useful, although elementary observations
on derivatives of $g_k$ and $x_k$ for $k\geq 0$.

The following formula is easy to check, for $a,b\in C$ not both vanishing:
\begin{equation}
{\cal D}_n((az+b)^{-w})=\binomial{w+n-1}{n}(-1)^n\frac{a^n}{(az+b)^{n+w}},
\label{eq:derivation_azb}\end{equation}
where the operators $\mathcal{D}_n=(-\widetilde{\pi})^nD_n$
have been introduced in \cite[Section 1]{BP}.

\begin{Lemma}\label{lemme:2ek} For $q\not=2,k\geq 1$ or for $q=2,k\geq 2$, we have: 
$$D_2g_k=\cdots=D_{q^k-1}g_k=0.$$
Moreover, for all $q$ and for all $k\geq 1$, we have:
$$D_1x_k=\cdots=D_{q^k-1}x_k=0.$$
\end{Lemma}
\noindent {\em Proof.} 
Assuming that $k\geq 2$ for $q=2$ or $k\geq 1$ otherwise,
the integer $w=q^k-1$ is $\geq 2$ and we have the following congruences:
\begin{small}
\[n=2,\qquad \binomial{w+n-1}{2}=\binomial{q^k}{2}\equiv0\pmod{p},\]
\[\begin{array}{cc}n=p^s,\\ s=1,\ldots,ek-1\end{array}\qquad
\binomial{w+n-1}{n}=\binomial{p^s+q^k-2}{p^s}
\equiv 0\pmod{p}.\]
\end{small}
Therefore, from the computations above, the uniform convergence of Eisenstein
series (\ref{eq:definition_eisenstein}) and  formula (\ref{eq:derivation_azb}):
\begin{equation}
{\cal D}_2E^{(q^{k}-1)}={\cal D}_{p}E^{(q^{k}-1)}={\cal D}_{p^2}E^{(q^{k}-1)}=
\cdots={\cal D}_{p^{ek-1}}E^{(q^{k}-1)}=0.\label{eq:derivaateE}\end{equation}
By (\ref{def_gk}) and the fact that $D$ is iterative (use formulas (18) and (19)
of \cite{BP}), we obtain the property concerning the $g_k$'s.

The property of the derivatives of the $x_k$'s then follows from
the definition $x_k=D_1g_k$. Indeed, $D_1x_k=D_1(D_1g_k)=2D_2g_k=0$
(this holds when $q=2$ by the congruence $2\equiv 0$ and when $q\not=2$
by the equality $D_2g_k=0$ we have just proved). 
Furthermore, for $2\leq i\leq q^k-1$, $D_ix_k=D_i(D_1g_k)=D_1(D_ig_k)=0.$ \CVD

We recall, in the next proposition and for the rest of the paper, that $L_k:=[k]\cdots[1]$ for $k>0$.
We also set $L_0:=1$.

\begin{Proposition}\label{xkordqk} For all $k\geq 0$, we have
\begin{equation}\label{beginningseriesxk}
x_k=(-1)^{k+1}L_ku^{q^k}+\cdots.
\end{equation}
In particular, $\nu_\infty(x_k)=q^k$.
\end{Proposition}
\noindent\emph{Proof.} By Lemma \ref{lemme:2ek} and \cite[Lemma~5.2]{BP} we have,
for all $k\geq 0$, $x_k\in C[[u^{q^k}]]$. Since 
$x_k$ vanishes at infinity, we may write:
$$x_k=\sum_{i\geq 1}c_{k,i}u^{iq^k},\quad c_{k,i}\in C,\quad k\geq 0.$$
From Corollaries (10.5), (10.11) and (10.4) of \cite{Ge}, collected
in the first part of Lemma 4.2 of \cite{BP}, we find the following $u$-expansions:
\begin{eqnarray*}
g&=&1-[1]u^{q-1}+\cdots\in C[[u^{q-1}]],\\
h&=&-u(1+u^{(q-1)^2}+\cdots)\in uC[[u^{q-1}]],\\
E&=&u(1+u^{(q-1)^2}+\cdots)\in u C[[u^{q-1}]].
\end{eqnarray*}
The third $u$-expansion tells that the result is true for $k=0$, since
$x_0=-E=-u+\cdots$. 
We also verify the result for $k=1$ because the three $u$-expansions yield
$x_1=-Eg-h=[1]u^q+\cdots$. 

We finish the proof by induction on $s=k-2\geq 0$
with the help of Proposition~\ref{sequences}.
From the recursion formula for $x_{s+2}$ we see 
that the coefficient of $u^{q^{s+1}}$ in the $u$-expansion of
$x_{s+2}$ is $c_{s+1,1}+[s+1]c_{s,1}$. But $x_{s+2}\in
C[[u^{q^{s+2}}]]$ and there cannot be a non-trivial contribution
by a term proportional to $u^{q^{s+1}}$; we deduce that this coefficient is zero.
Therefore, $c_{s+1,1}=-[s+1]c_{s,1}$.\CVD

\begin{Remark}{\em It can be proved that for all $k$, the normalisation of $x_k$ lies in $A[[u^{q^k}]]$.}\end{Remark}

\subsubsection{Tables.\label{tables}}\label{sectiontables}

The following table collects several useful data checked above. The
order of vanishing of $h_k$ easily follows from (\ref{eq:egkhk}), the fact that
$\nu_\infty(g_k)=0$ and Proposition~\ref{xkordqk}.
The index $k$ is supposed to be $\geq 1$.
\begin{center}
\medskip
\begin{tabular}{|l|l|l|l|l|}
\hline
form $f$ & weight & depth & type & $\nu_\infty(f)$ \\
\hline
$g_k$ & $q^k-1$ & $0$ & $0$ & $0$ \\
\hline
$h_k$ & $q^k+1$ & $0$ & $1$ & $1$ \\
\hline
$x_k$ & $q^k+1$ & $1$ & $1$ & $q^k$ \\
\hline
\end{tabular}
\end{center}
\medskip
The next table describes the first values of $x_k$, from which one easily deduces
the corresponding values of $g_k,h_k$ thanks to (\ref{eq:egkhk}):
\begin{small}
\begin{center}
\medskip
\begin{tabular}{|rcl|}
\hline
$x_0$ & $=$ & $-E$ \\
\hline
$x_1$ & $=$ & $-Eg-h$\\
\hline
 $x_2$ & $=$ & $-E(g^{q+1}+[1]h^{q-1})-g^q h$\\
\hline
\end{tabular}
\medskip
\end{center}
\end{small}

\subsubsection{A multiplicity estimate for forms of depth $\leq q$}\label{sectionmultxk}
 The next simple lemma will be used quite often.

\begin{Lemma}\label{lemme_resultant}
Let $f,g$ be quasi-modular forms, with $f\in\widetilde{M}^{\leq l}_{w,m}$ and $g\in\widetilde{M}^{\leq l'}_{w',m'}$, considered 
as polynomials in $C[E,g,h]$.
Then, their resultant $\rho:=\mathbf{Res}_E(f,g)$ with respect to $E$ is a Drinfeld modular form of weight
$w(\rho)=lw'+wl'-2ll'$ and type $m(\rho)=lm'+l'm-ll'$.
\end{Lemma}

\noindent\emph{Proof.} This is elementary and follows by a suitable adaptation of, for example, \cite[Lemma 6.1]{PF}
(see also \cite[Theorem 6.1]{Sturmfels}). The information on the type will not be used in this paper but is given for the sake of completeness.\CVD

The {\em degree} $d(f)$ of a quasi-modular form $f$ is by definition the positive integer $d(f):=w(f)-l(f)$,
difference between its weight and its depth. 

As an application of the previous results, we prove here a
multiplicity estimate for quasi-modular forms of depth $\leq q$
that will be used later in this paper.

\begin{Proposition}\label{th_multiplicity} Let $w$ and $m$ be integers
such that $\widetilde{M}^{\leq q}_{w,m}\not=\{0\}$, and let $f$ be a
non-vanishing form in $\widetilde{M}^{\leq q}_{w,m}$. Then
$$\nu_\infty(f)\leq \frac{q^2+1}{q+1}d(f).$$
\end{Proposition}

\noindent\emph{Proof.} If the bound of the Proposition holds for
two forms $f_1$ and $f_2$ with $f_1f_2\in\widetilde{M}^{\leq q}_{w,m}$, then
the bound clearly holds for $f_1f_2$ too, by adding the inequalities.
So we may suppose that $f$ is irreducible in the ring $C[E,g,h]$.
Let $k$ be the smallest integer
$\ge 0$ such that
\begin{equation}\label{choix_de_k}
w(f)<q^k+q.
\end{equation}
If there is $\lambda\in C$ such that $f=\lambda x_k$,
then $\nu_{\infty}(f)=q^k$ and $d(f)=q^k$, so the bound of
Proposition~\ref{th_multiplicity} holds. If it is not the case,
then consider the resultant (with respect to the indeterminate $E$)
$$\rho:=\mathbf{Res}_E(f,x_k).$$
The function $\rho$ is a non zero modular form of weight
$w(f)+l(f)(q^k-1)$ according to Lemma \ref{lemme_resultant}, so we have,
by \cite[formula (5.14)]{Ge},
\begin{equation}\label{majo}
\nu_{\infty}(\rho)\le w(\rho)/(q+1)\le (w(f)+q(q^k-1))/(q+1).
\end{equation}
On the other hand, since there exist $U, V\in\widetilde{M}$ such that
$\rho=Uf+Vx_k$, we have
\begin{equation}\label{mino}
\nu_{\infty}(\rho)\ge
\min\{\nu_{\infty}(f),\nu_{\infty}(x_k)\} = \min\{\nu_{\infty}(f),q^k\}.
\end{equation}
By (\ref{choix_de_k}) we have $q^k>(w(f)+q(q^k-1))/(q+1)$, so the
compatibility of (\ref{majo}) and (\ref{mino}) implies
$\min\{\nu_{\infty}(f),q^k\}=\nu_{\infty}(f)$, hence
\begin{equation}\label{maj}
\nu_{\infty}(f)\le (w(f)+q(q^k-1))/(q+1).
\end{equation}
We distinguish now two cases. If $k=0$, then we get (note that $d(f)>0$)
$$(q+1)\nu_{\infty}(f)\le w(f)\le d(f)+q\le (q+1) d(f)$$
and the result follows in this case. If $k\ge 1$, then by minimality
of $k$ satisfying (\ref{choix_de_k}) we have $q^{k-1}+q\le w(f)$,
hence $q^k\le q(w(f)-q)\le q d(f)$. Replacing in (\ref{maj}) and
using the estimate $w(f)\le d(f)+q$, we get the result. \CVD

\begin{Remark}\label{remark_multiplicity}{\em  Observing the end of the proof of this Proposition,
we notice that if $l(f)<q$, then we have the strict inequality
$w(f)<d(f)+q$, and hence we get the strict inequality
$\nu_{\infty}(f)<(q^2+1)/(q+1)\, d(f)$. This remark will be crucial in the proof
of Lemma~\ref{xik_irreducible}. }\end{Remark}

\subsubsection{Proof of Theorem \ref{theo_xk_extremal}.}\label{extremality_xy}

The first table of Section~\ref{sectiontables} shows that $x_k$ has
a high vanishing order compared to its weight. Here we prove that
this vanishing order is the highest possible among forms in
$\widetilde{M}^{\leq 1}_{q^k+1,1}$, that is,
that the form $x_k$ is extremal in $\widetilde{M}^{\leq 1}_{q^k+1,1}$.
In fact we will even get a slightly more general result
(Proposition~\ref{x0x1x2}).

\begin{Lemma}\label{lemma:resultants}
For all $k\geq 0$ we have:
\[\rho_{k}:=\det\sqm{g_k}{h_k}{g_{k+1}}{h_{k+1}}=(-1)^kL_kh^{q^k}.\]
\end{Lemma}


\noindent {\em Proof.} 
To compute $\rho_{k}$ it suffices to substitute $g_{k+1},h_{k+1}$ by their
expressions as linear combinations of $g_k,g_{k-1}$ and $h_k,h_{k-1}$ with
coefficients in $M$ (\textit{cf.} (\ref{eq:g0g1gk}) and
Proposition~\ref{sequences}). We obtain the formula
$\rho_{k}=-[k]h^{q^{k-1}(q-1)}\rho_{k-1}$; since $\rho_{0}=
\det\sqm{1}{0}{g}{h}$, the lemma follows.

\CVD

%

\begin{Lemma}\label{xkirreducible}
For all $k\geq 0$, the form $x_k$ is irreducible as a polynomial of $C[E,g,h]$.
\end{Lemma}

\noindent\emph{Proof.} 
Assume by contradiction that $\delta$ is a non-trivial divisor 
of the polynomial $x_k$. Since $x_k=-g_kE-h_k$ is of depth $1$, we can assume without loss of generality that
$\delta$ is a modular form,
common divisor of $g_k$ and $h_k$. But then, $\delta$ divides the form
$\rho_{k}$ of Lemma \ref{lemma:resultants}, which
tells that $\delta$ is a multiple by an element of $C^\times$ of a power of 
$h$. Hence, $h$ divides $g_k$, which does not vanish at infinity;
contradiction, because $h$ does. \CVD

The next Lemma gives a sufficient condition for an irreducible quasi-modular form
to be extremal.

\begin{Lemma}\label{lemmaextremal}
Let $\varphi\in\widetilde{M}_{w,m}^{\leq l}$ be a quasi-modular form satisfying
\begin{equation}\label{conditionextremal}
(q+1)\nu_{\infty}(\varphi)>2 l(\varphi)d(\varphi).
\end{equation}
Then, for all non zero quasi-modular forms $f\in\widetilde{M}_{w,m}^{\leq l(\varphi)}$
without non-constant common factor with $\varphi$, we have
$\nu_{\infty}(f)\leq \nu_{\infty}(\varphi)$.
In particular, if $\varphi$ is irreducible in $C[E,g,h]$ then it is extremal.
\end{Lemma}

\noindent\emph{Proof.} Suppose that there exists a non zero form
$f\in \widetilde{M}^{\le l(\varphi)}_{w,m}$
such that $\nu_{\infty}(f)>\nu_{\infty}(\varphi)$, and
such that $f$ and $\varphi$ have no common factor.
Then the resultant $\rho:=\mathbf{Res}_E(f,\varphi)$ is non zero.
Note that $\rho$ is a modular form of weight
$w(\rho)=wl(\varphi)+wl(f)-2l(\varphi)l(f)\le 2l(\varphi)d(f)$ (Lemma \ref{lemme_resultant}).
Now, we have on one side $\nu_{\infty}(\rho)\le w(\rho)/(q+1)$
(since $\rho$ is modular), and on the other side $\nu_{\infty}(\rho)\ge
\min\{\nu_{\infty}(\varphi),\nu_{\infty}(f)\}=\nu_{\infty}(\varphi)$. Thus we find
$$
\nu_{\infty}(\varphi)\le
w(\rho)/(q+1)\le 2l(\varphi)d(f)/(q+1).
$$
But this contradicts the hypothesis (\ref{conditionextremal}). This shows that
$\nu_{\infty}(f)\leq \nu_{\infty}(\varphi)$ and the first part of the
Lemma is proved. The second one is clear, since when $\varphi$ is irreducible,
then any $f\in\widetilde{M}^{\le l(\varphi)}_{w,m}$ either has no common
factor with $\varphi$ or has the form $\lambda \varphi$ with $\lambda\in C$. \CVD

\noindent {\em Proof of Theorem \ref{theo_xk_extremal}.}
It follows at once from
Lemma~\ref{lemmaextremal} applied with $\varphi=x_k$, since $x_k$ is irreducible by
Lemma~\ref{xkirreducible} and since the condition (\ref{conditionextremal})
is clearly satisfied (see the first table of Section~\ref{sectiontables}). \CVD

\begin{Remark}{\em 
We can generalise Theorem \ref{theo_xk_extremal} a little bit as follows.

\begin{Proposition}\label{x0x1x2} Let $r_0,\ldots,r_s$ be non-negative integers
not all of which are zero, let us write
$l=r_0+\cdots+r_s,w=r_0+r_1q+\cdots+r_sq^s+l$ and let $m$ be the class
of reduction modulo $q-1$ of $l$.
If $l<\frac{q+1}{2}$, then the quasi-modular form 
$$x=x_0^{r_0}\cdots x_s^{r_s}\in\widetilde{M}^{\leq l}_{w,m}
\setminus\widetilde{M}^{\leq l-1}_{w,m}$$ is extremal.
\end{Proposition}

\noindent\emph{Proof.} 
Let $f$ be a non-zero element of $\widetilde{M}^{\leq l}_{w,m}$.
Write $f=\varphi^*\delta$ and $x=\varphi\delta$, where $\delta$, $\varphi$, $\varphi^*$
are elements of $C[E,g,h]$ such that $\varphi$ and $\varphi^*$ are coprime.
If $\varphi$ is constant, then $f$ is a multiple of $x$ in $C[E,g,h]$, so
$f=\lambda x$ for some $\lambda\in C^*$ (since $w(f)=w(x)=w$ by hypothesis).
Thus $\nu_\infty(f)=\nu_\infty(x)$ in this case.
If now $\varphi$ is not constant, then $\varphi$ has the form (up to an
element of $C^*$) $$\varphi=x_0^{\alpha_0}\cdots x_s^{\alpha_s}$$
with $\alpha_0,\ldots,\alpha_s$ not all zero such that
$0\leq\alpha_i\leq r_i$ for all $i$.
One readily checks that the condition
$(q+1)\nu_{\infty}(\varphi)>2 l(\varphi)d(\varphi)$
of Lemma~\ref{lemmaextremal} is equivalent to  $q+1>2l$, so it is
satisfied. Applying now this lemma yields
$\nu_\infty(\varphi^*)\leq\nu_\infty(\varphi)$, or equivalently
$\nu_\infty(f)\leq\nu_\infty(x)$. Thus $x$ is extremal. \CVD}\end{Remark}


\subsection{The family $(\xi_k)_{k\geq 0}$.\label{sectionxik}}

We recall that the definition of the forms $\xi_k$ occurs in (\ref{def_xi_k}).
In this section we study the forms $\xi_k$ and give proofs of Theorems \ref{xik_extremal} and \ref{th_multiplicity3}.


\begin{Proposition}\label{wlnuxik}
For all $k\geq 0$, the quasi-modular form $\xi_k$ satisfies
$$
w(\xi_k)=(q^k+1)(q+1),\quad l(\xi_k)=q+1,\quad \nu_{\infty}(\xi_k)=q^{k+2}+q^k.
$$
\end{Proposition}

\noindent\emph{Proof.}
The fact that $w(\xi_k)=(q^k+1)(q+1)$ immediately
follows from the definition of $\xi_k$. Let us prove that the depth of $\xi_k$
is $q+1$. If $k=0$, a straightforward computation yields
$\xi_0=-[1]E^{q+1}+ghE+h^2$, so the result is clear. Suppose now that
$k\ge 1$. By (\ref{eq:egkhk}), $x_k=-g_kE-h_k$ for all $k$.
We get, by definition of $\xi_k$:
\begin{equation}\label{xxxik}
\xi_k= [k]^q(g_{k-1}^qE^q+h_{k-1}^q)(g_{k+1}E+h_{k+1})
- [k+1](g_kE+h_k)(g_k^qE^q+h_k^q).
\end{equation}
If we consider $\xi_k$ as a polynomial in $E$, the coefficient of $E^{q+1}$
is therefore equal to $\alpha:=[k]^qg_{k+1}g_{k-1}^q-[k+1]g_k^{q+1}$.
This last form is non zero
since the constant term of its $u$-expansion is (using the fact that
$g_i=1+\cdots$ for all $i$) $[k]^q-[k+1]=-[1]\not=0$.
It follows that the degree
in $E$ of the form $\xi_k$ is exactly $q+1$, hence $l(\xi_k)=q+1$.

It remains to prove that $\nu_{\infty}(\xi_k)=q^{k+2}+q^k$.
To do this, we first notice that the following relation holds :
\begin{equation*}
-\Delta^{q^k}\xi_k=x_{k+1}^{q+1}-x_k^qx_{k+2}.
\end{equation*}
Indeed, using the recursion formula of the sequence $(x_k)$
(Proposition~\ref{sequences}), we have:
\begin{eqnarray*}
x_{k+1}^{q+1}-x_k^qx_{k+2} & = & (g^{q^{k+1}}x_k^q-
[k]^q\Delta^{q^k}x_{k-1}^q)x_{k+1}\cr 
& & \phantom{g^{q^{k+1}}x_k^q} -x_k^q(g^{q^{k+1}}x_{k+1}-[k+1]
\Delta^{q^k}x_k)\cr
& = & -\Delta^{q^k}\xi_k.
\end{eqnarray*}
Thus, it suffices to show that
\begin{equation}\label{ordxik}
\nu_{\infty}(x_{k+1}^{q+1}-x_k^qx_{k+2})=q^{k+2}+q^{k+1}.
\end{equation}
But by Proposition~\ref{xkordqk} we have:
\begin{eqnarray*}
x_{k+1}^{q+1}-x_k^qx_{k+2}
& = & \bigl((-1)^{k+2}L_{k+1}u^{q^{k+1}}+\cdots \bigr)^{q+1}\cr
& & \ {}- \bigl((-1)^{k+1} (-1)^{k+3}L_k^qL_{k+2}
u^{q^{k+1}+q^{k+2}} + \cdots\bigr)\cr
& = & ([k+1]^q-[k+2])[k+1]L_{k}^{q+1} u^{q^{k+1}+q^{k+2}}
+\cdots\cr
& = & -[k+1][1]L_k^{q+1} u^{q^{k+1}+q^{k+2}}+\cdots
\end{eqnarray*}
(we have used the fact that $[k+1]^q-[k+2]=-[1]$).
Hence (\ref{ordxik}) holds and the Proposition is proved.
\CVD

\subsubsection{Proof of Theorem \ref{xik_extremal}.}

\begin{Lemma}\label{xik_irreducible}
For every $k\ge 0$ the form $\xi_k$ is irreducible in $C[E,g,h]$.
\end{Lemma}

\noindent\emph{Proof.} Suppose that $\xi_k$ is reducible. Write
$\xi_k=ab$, where $a$, $b$ are non constant quasi-modular forms.
Suppose first that $l(a)\ge 1$ and $l(b)\ge 1$, or, equivalently, that
$l(a)\le q$ and $l(b)\le q$. Since $l(\xi_k)= q+1$,
$a$ and $b$ cannot both have a depth equal to $q$, so
we may suppose $l(a)<q$. By Proposition~\ref{th_multiplicity}
and Remark \ref{remark_multiplicity}, we have
$$ \nu_{\infty}(a)< \frac{q^2+1}{q+1}\, d(a)\quad {\rm and}
\quad \nu_{\infty}(b)\le \frac{q^2+1}{q+1}\, d(b). $$
Hence, since $d(a)+d(b)=d(\xi_k)= q^k(q+1)$,
$$ \nu_{\infty}(\xi_k)= \nu_{\infty}(a)+\nu_{\infty}(b)<\frac{q^2+1}{q+1}\,
(d(a)+d(b))=q^k(q^2+1).$$
But this contradicts the fact that $\nu_{\infty}(\xi_k)=q^k(q^2+1)$
(Proposition~\ref{wlnuxik}).

Thus we have $l(a)=0$ or $l(b)=0$. We will suppose in what follows
that $l(a)=0$, \ie\ $a$ is a modular form. We will even assume,
without loss of generality, that $a$ is irreducible. Returning to the expression
(\ref{xxxik}) of $\xi_k$, we see that $\xi_k=\alpha E^{q+1}+\beta E^q+\gamma E+\delta$ with
\begin{eqnarray*}
\alpha&=&[k]^q g_{k-1}^qg _{k+1} - [k+1]g_k^{q+1},\cr
\beta&=&[k]^q g_{k-1}^qh _{k+1} - [k+1]g_k^qh_k,\cr
\gamma&=&[k]^q h_{k-1}^qg _{k+1} - [k+1]h_k^qg_k,\cr
\delta&=&[k]^q h_{k-1}^qh _{k+1} - [k+1]h_k^{q+1}.
\end{eqnarray*}
As in Lemma~\ref{lemma:resultants} we define
$$\rho_{k}:=\det\sqm{g_k}{h_k}{g_{k+1}}{h_{k+1}}.$$
We have (Lemma~\ref{lemma:resultants}) $\rho_k=(-1)^kL_kh^{q^k}$.
Since $a$ is a modular form dividing $\xi_k$, it divides $\alpha$, $\beta$, $\gamma$ and $\delta$,
and thus also the two forms
$$ h_{k-1}^q\alpha-g_{k-1}^q\gamma= [k+1]g_k\rho_{k-1}^q= (-1)^{k-1}L_{k-1}^q
[k+1] h^{q^{k}}g_k$$
and
$$ h_{k-1}^q\beta-g_{k-1}^q\delta= [k+1]h_k\rho_{k-1}^q= (-1)^{k-1}L_{k-1}^q
[k+1] h^{q^{k}}h_k.$$
But $h$ does not divide $\alpha$ as this form does not vanish at infinity (see \eg\
the proof of Proposition~\ref{wlnuxik}). So $a$ must divide both
$g_k$ and $h_k$. But then $a$ divides the form $\rho_k$, hence is equal to
$h$ (up to a constant factor). But this is impossible since $h$ does not divide $\alpha$.
Finally, the contradiction obtained shows that $\xi_k$ is irreducible,
as announced.\CVD

\noindent\emph{Proof of Theorem \ref{xik_extremal}.} By Proposition~\ref{wlnuxik} we have
$\nu_{\infty}(\xi_k)=q^k(q^2+1)$, $l(\xi_k)=q+1$ and $d(\xi_k)=q^k(q+1)$,
hence the condition $(q+1)\nu_{\infty}(\xi_k)>2l(\xi_k)d(\xi_k)$
is satisfied for $q>2$. Since $\xi_k$ is irreducible by
Lemma~\ref{xik_irreducible},
we may therefore apply Lemma~\ref{lemmaextremal}.
We get that $\xi_k$ is extremal.\CVD

\begin{Remark}{\em For $q=2$, numerical computations show that $\xi_k$ is also extremal for $k=0,1$ (cf. Section \ref{experiments}).}\end{Remark}

\subsubsection{Proof of Theorem \ref{th_multiplicity3}.}\label{sectionmultxik}

We argue as in the proof of
Proposition~\ref{th_multiplicity}.
We may suppose that $f$ is irreducible in the ring $C[E,g,h]$.
Let $k$ be the smallest integer
$\ge 0$ such that
\begin{equation}\label{choix_de_k2}
w(f)<q^k+q^2.
\end{equation}
If there is $\lambda\in C$ such that $f=\lambda \xi_k$,
then $\nu_{\infty}(f)=q^k(q^2+1)$ and $d(f)=q^k(q+1)$, so the bound of
the theorem holds. If it is not the case,
then consider the resultant $\rho:=\mathbf{Res}_E(f,\xi_k)$.
The function $\rho$ is a non zero modular form of weight
$(w(f)+l(f)(q^k-1))(q+1)$ by Lemma \ref{lemme_resultant}, so we have
\begin{equation}\label{majo2}
\nu_{\infty}(\rho)\le w(\rho)/(q+1)\le w(f)+q^2(q^k-1).
\end{equation}
On the other hand, we have
\begin{equation}\label{mino2}
\nu_{\infty}(\rho)\ge
\min\{\nu_{\infty}(f),\nu_{\infty}(\xi_k)\} = \min\{\nu_{\infty}(f),q^k(q^2+1)\}.
\end{equation}
By (\ref{choix_de_k2}) we have $q^k(q^2+1)>w(f)+q^2(q^k-1)$, so the
compatibility of (\ref{majo2}) and (\ref{mino2}) implies
$\min\{\nu_{\infty}(f),q^k(q^2+1)\}=\nu_{\infty}(f)$, hence
\begin{equation}\label{maj2}
\nu_{\infty}(f)\le w(f)+q^2(q^k-1).
\end{equation}
We now distinguish two cases. If $k=0$, then we get (notice that $d(f)>0$)
$$\nu_{\infty}(f)\le w(f)\le d(f)+q^2\le (q^2+1) d(f)$$
and the result follows in this case. If $k\ge 1$, then by minimality
of $k$ satisfying (\ref{choix_de_k2}) we have $q^{k-1}+q^2\le w(f)$,
hence $q^k\le q(w(f)-q^2)\le q d(f)$. Replacing in (\ref{maj2}) and
using the estimate $w(f)\le d(f)+q^2$, we get the result.\CVD

\section{Differential extremality.\label{differential_extremality}}

In Section 5 of \cite{BP}, we have introduced the following
subgroups of $F^\times$ where $F=C(E,g,h)$:
$$
\Psi_k =\{f\in F^\times;\ (D_{p^j}f)/f\in\widetilde{M}
\mbox{ for all }
0\leq j\leq k\},$$
with the additional notation $\Psi_{-1}:=F^\times$. The differential exponent (the map we have mentioned 
in the introduction) can then be defined in the following alternative way:
$$\epsilon_D:F^\times\rightarrow\ZZ_{\geq 0}\cup\{\infty\}$$
such that $\epsilon_D(f)=k+1$ if $f\in\Psi_{k}\setminus\Psi_{k+1}$ and $\epsilon_D(f)=\infty$
if $f\in\Psi_\infty:=\cap_{i=-1}^\infty\Psi_i$. 
This makes sense because for all $k\geq -1$, $\Psi_k\supsetneq\Psi_{k+1}$.

\begin{Lemma}\label{lemma:epsilonD}
The following properties of the differential exponent hold, with $f,g\in F^\times$.
\begin{enumerate}\item If $f\in C^\times h^\ZZ$, then $\epsilon_D(f)=\infty$.
\item If $q\not=2,3$ then $\epsilon_D(f)=\infty$ implies that $f\in C^\times h^\ZZ$.
\item We have $\epsilon_D(f^p)=\epsilon_D(f)+1$. Moreover, if $p\nmid m$, $\epsilon_D(f^m)=\epsilon_D(f)$. 
\item We have $\epsilon_D(fg)\geq \inf\{\epsilon_D(f),\epsilon_D(g)\}$, and equality holds when $\epsilon_D(f)\neq \epsilon_D(g)$.
\end{enumerate}
\end{Lemma}
\emph{Sketch of proof.} The first property follows from \cite[Proposition 3.6]{BP}. The second property is a paraphrase of 
\cite[Theorem 3]{BP}. In the third property, the first part is clear by
using the formula $D_{p^{k+s}}f^{p^k}=(D_{p^s}f)^{p^k}$, which follows easily from Leibniz rule \cite[Equation (16)]{BP}.
As for the second part,
we observe that $\epsilon_D(f^m)\geq \epsilon_D(f)$ because the sets $\Psi_k$
are, as already observed, multiplicative subgroups of $F^\times$.
Let us write $k=\epsilon_D(f)$. Then, by Leibniz's formula (15) of \cite{BP},

$$\frac{D_{p^k}(f^m)}{f^m}=m\frac{D_{p^k}f}{f}+
\sum_{\genfrac{}{}{0pt}{2}{i_1+\cdots+i_m=p^k}{0\leq i_1,\ldots,i_m<p^k}}
\frac{D_{i_1}f}{f}\cdots\frac{D_{i_m}f}{f}.$$
The sum on the right hand side is an element of $\widetilde{M}$
while $(D_{p^k}f)/f\not\in\widetilde{M}$, which
implies the inequality $\epsilon_D(f^m)\leq \epsilon_D(f)$.
The fourth property can be proved in a similar way; its proof is then
left to the reader.\CVD

As in \cite[Section~5]{BP} we denote by $F_r$ the following subset
of $F$, which turns out to be a subfield:
$$ F_r=\bigcap_{0\le i\le r}\ker D_{p^i}.$$
For $r=-1$ we define $F_{-1}=F$. We have $F_{r+1}\subset F_r$ and
$F_r^p\subset F_{r+1}$ for all $r\ge 0$.

\begin{Proposition}\label{f0hs} 
Let $f$ be an element of $\Psi_k$ ($k\geq 0$). Then, $h^{-\nu_\infty(f)}f\in F^\times_k$. More precisely,
$$\frac{D_nf}{f}=\frac{D_n(h^{\nu_\infty(f)})}{h^{\nu_\infty(f)}},\quad n=0,\ldots,p^{k+1}-1.$$
\end{Proposition}

\medskip

\noindent\emph{Proof}.  Let us write $\nu_\infty(f)=n_0+\cdots+n_kp^k+mp^{k+1}$ with 
$n_0,\ldots,n_k\in\{0,\ldots,p-1\}$ and $m\geq 0$. Let us also define inductively:
$$f_{-1}=f,f_0=f_{-1}/h^{n_0},\ldots,f_{s}=f_{s-1}/h^{n_sp^s},\ldots$$
By \cite[Proposition 3.6]{BP}, we know that $f_{-1},f_0,\ldots,f_k\in\Psi_k$. 

We now prove by induction on $s$ that $f_s\in F_s^{\times}$ for
$0\leq s\leq k$.
From elementary weight considerations, there exists $\alpha\in C$ such that 
$$D_1f_{-1}=\alpha Ef_{-1}.$$
Writing $f_{-1}=cu^{i}+\cdots$ (the dots $\ldots$ are understood as a series of higher powers of $u$ and $c\neq 0$)
and comparing $D_1f_{-1}=ciu^{i+1}+\cdots$ with $\alpha Ef_{-1}=\alpha
(u+\cdots)(cu^{i}+\cdots)$
yields $\alpha \equiv i\equiv n_0\pmod{p}$. Hence $\alpha\in\FF_p=\ZZ/p\ZZ$
and we can choose a representative $\alpha$
of this class, with $\alpha=n_0$ (allowing an abuse of notation).
Now,
$$D_1h^{n_0}=n_0Eh^{n_0},$$
which implies $D_1f_0=0$, hence $f_0\in F_0$.

Assuming now that $f_s\in F^\times_s$ for $s<k$,
we proceed to prove that $f_{s+1}\in F_{s+1}^\times$. Since $f_s\in\Psi_s$,
we can apply Lemma 5.9 of \cite{BP} to 
check that:
$$D_{p^{s+1}}f_s=\alpha_s E^{p^{s+1}}f_s,\quad \alpha_s\in\{0,\ldots,p-1\}.$$
Since $f_s\in F_s^\times$ by hypothesis, Lemma 5.2 of \cite{BP} implies that
$$f_s=c_su^{i_sp^{s+1}}+\cdots\in C((u^{p^{s+1}})),\quad c_s\neq 0.$$
Comparing $D_{p^{s+1}}f_s=c_si_su^{(i_s+1)p^{s+1}}+\cdots$ with
$\alpha_sE^{p^{s+1}}f_s=\alpha_s(u^{p^{s+1}}+\cdots)(c_su^{i_sp^{s+1}}+\cdots)$
yields $\alpha_s\equiv i_s\equiv n_{s+1}\pmod p$ which implies $\alpha_s=n_{s+1}$.
Now, 
$$D_{p^{s+1}}h^{n_{s+1}p^{s+1}}=n_{s+1}E^{p^{s+1}}h^{n_{s+1}p^{s+1}}$$
so that $f_{s+1}\in F^\times_{s+1}$.

Finally,
$$\frac{f}{h^{\nu_\infty(f)}}=\frac{f}{h^{n_0+n_1p+\cdots+n_kp^k}}\frac{1}{h^{mp^{k+1}}}=\frac{f_k}{h^{mp^{k+1}}}\in F_k^\times,$$
and for all $0\leq n\leq p^{k+1}-1$, $$\frac{D_{n}f}{f}=\frac{D_n(h^{\nu_\infty(f)})}{h^{\nu_\infty(f)}}.$$
\CVD

\subsection{Differential extremality of the forms $x_k$}

Here we prove Theorem \ref{theorem:xkbis}.
The proof makes use at once of the main result of \cite{GG}, where the 
factorisations in $K[g,h]$ of Gekeler's forms $g_k$'s are determined,
and of differential tools that will be
developed in Section~\ref{differential_tools}. We begin with a 
subsection devoted to the description of the needed result from \cite{GG},
and to some preliminary lemmas.

\subsubsection{Preliminaries to the proof of Theorem \ref{theorem:xkbis}}

Since the $C$-algebra $\widetilde{M}$ is equal to the
polynomial ring $C[E,g,h]$,
the group $\mathbf{Aut}(C/K)$ naturally acts on $\widetilde{M}$,
the action on an element $f\in\widetilde{M}$ being
given by the action on the coefficients of $f$ seen as a polynomial in $E,g,h$.
Since the functions $E,g,h$ have their $u$-expansions with coefficients in $K$,
if $\sigma$ is a $K$-automorphism and if $f=\sum_{i\geq 0}c_iu^i$, then
we have
$$f^\sigma=\sum_{i\geq 0}c_i^\sigma u^i.$$
If $L$ is a subfield of $C$, we will say that
$f\in\widetilde{M}$ is {\em defined over $L$} if $f\in L[E,g,h]$.

\begin{Lemma}\label{inseparable}
For $k\geq 0$, let $f$ be normalised and differentially extremal in $\widetilde{M}^{\leq 1}_{q^k+1,1}$.
Then, $f$ is defined over the inseparable closure of $K$ in $C$.
\end{Lemma}

\noindent\emph{Proof.}
Let $f$ be as in the statement of the lemma.
We know that $\epsilon_D(f)\geq e(k+1)=\epsilon_D(x_k)$
by Proposition \ref{propformulasgsxs}. Hence,
for all $n$ such that $1\leq n\leq q^{k+1}-1$, $f$ divides $D_nf$.
Let $A_n\in\widetilde{M}$ be the polynomial such that
$D_nf=A_n f$, for $n$ as above.
By Proposition \ref{f0hs}, there exists $\nu\in\ZZ$
such that $A_n=(D_nh^\nu)/h^{\nu}$, so $A_n$ is defined over $K$ by
\cite[formula (3)]{BP}.
Let $\sigma$ be any  element of $\mathbf{Aut}(C/K)$. For 
all $n$ with $1\leq n\leq q^{k+1}-1$, we have:
\begin{equation*}
D_nf^\sigma=(D_nf)^\sigma=(A_nf)^\sigma =A_nf^\sigma,
\end{equation*}
where the first identity comes again from \cite[formula (3)]{BP}.
Define $\nu:=\nu_\infty(f)=\nu_\infty(f^\sigma)$ and suppose that
$f\not=f^{\sigma}$.
By Proposition \ref{f0hs} and \cite[Lemma 5.2]{BP}, we
have $fh^{-\nu}, f^{\sigma}h^{-\nu}\in C((u^{q^{k+1}}))$, so
there exists $\alpha\in\ZZ$ such that
$\nu_\infty((f-f^\sigma)h^{-\nu})=\alpha q^{k+1}$.
Now, since $f$ and $f^\sigma$ are normalised,
the form $\phi:=f-f^\sigma$ satisfies $\nu_\infty(\phi)>\nu$. Hence
we get $\alpha>0$, which implies $\nu_\infty(\phi)=\nu+\alpha q^{k+1}>q^k$.
But since
$\nu_{\infty}(x_k)=q^k$, this contradicts the extremality of
$x_k$ (Theorem~\ref{theo_xk_extremal}). Hence we have $f=f^\sigma$. \CVD

The following conjecture seems plausible, although pretty difficult to 
reach with the tools at our disposal.

\begin{Conjecture}
If $f$ is differentially extremal, then $f$ is defined over $K$.
\end{Conjecture}

We rewrite the needed result from \cite{GG} in the following lemma.

\begin{Lemma}\label{propocornelissen}
Let $k\geq 2$ be even. Then the polynomial
$g_k\in K[g,h]$ is irreducible and totally decomposable over a separable
extension of $K$.
\end{Lemma}
\noindent\emph{Proof.} 
We know from \cite[Section 2]{GG} that if we set $j=g^{q+1}/\Delta$,
there exist a non-zero element $c_k\in K$ and a polynomial $P_k\in K[X]$
such that
\begin{equation}
g_k=c_k\Delta^{\deg(P_k)}P_k(j),\label{eq:Pkj}
\end{equation}
where the degree of $P_k$ is
$$\frac{q^k-1}{q^2-1}.$$
Hence $p\nmid\deg(P_k)$, so that $P_k$ is separable over $K$.
Now, formula (\ref{eq:Pkj}) establishes a bijective
correspondence between the irreducible factors in $\overline{K}[X]$
of $P_k$ ($\overline{K}$ being an algebraic closure of $K$ in $C$), and
the irreducible factors in $\overline{K}[g,h]$ of $g_k$.
This shows the second statement of the lemma.
The statement on the irreducibility of $g_k$ follows
from \cite[Theorem E2]{GG}.\CVD

\begin{Remark}{\em If $k\geq 3$ is odd, it is easy to show from
the results in \cite{GG} that $g_k$ is
totally decomposable over a separable extension of $K$, and
equal to a product of $g$ and an irreducible polynomial of $K[g,h]$.
However, we will only use the case $k$ even.
}\end{Remark}

\begin{Lemma}\label{lemmaresultant}
Let $\mathcal{K}$ be any field. Let
$a,b,c,d$ be polynomials in $\mathcal{K}[X_1,\ldots,X_n]$, where $n\geq 1$,
and write $F=aX_0+b,G=cX_0+d$.
Let $\mathcal{I}=(F,G)$ be the ideal generated by $F,G$ in
$\mathcal{K}[X_0,X_1,\ldots,X_n]$, and define
$\mathcal{J}:=\mathcal{I}\cap\mathcal{K}[X_1,\ldots,X_n]$.
Then, if $a,c$ are coprime in $\mathcal{K}[X_1,\ldots,X_n]$,
the ideal $\mathcal{J}\subset\mathcal{K}[X_1,\ldots,X_n]$ is
principal, generated by the resultant $\rho:=\mathbf{Res}_{X_0}(F,G)=ad-bc$.
\end{Lemma}

\noindent\emph{Proof.} Let $AF+BG\in\mathcal{K}[X_1,\ldots,X_n]$ be
any polynomial of $\mathcal{J}$, where $A,B \in\mathcal{K}[X_0,\ldots,X_n]$.
Let us write
$$A=\sum_{i=0}^m\alpha_iX_0^i,\quad B=\sum_{i=0}^{m'}\beta_iX_0^i$$
with $\alpha_i,\beta_j\in \mathcal{K}[X_1,\ldots,X_n]$ and
$\alpha_m\beta_{m'}\neq 0$. Considering
the leading term of the expansion of $AF+BG$ in powers of $X_0$,
we see that $m=m'$. Moreover, assuming $m>0$ we have
\begin{equation}\label{c1}
\alpha_ma+\beta_mc=0
\end{equation}
and
\begin{equation}\label{c2}
\alpha_i b+\alpha_{i-1}a +\beta_i d + \beta_{i-1}c = 0,\quad 1\le i\leq m.
\end{equation}

Since $a,c$ are coprime, we deduce from (\ref{c1}) that there exists
$u_m\in \mathcal{K}[X_1,\ldots,X_n]$ such that 
$\alpha_m=cu_m$, $\beta_m=-au_m$. Substituting in
(\ref{c2}) for $i=m$, we get similarly the existence of
$u_{m-1}\in \mathcal{K}[X_1,\ldots,X_n]$ such that
$$
\alpha_{m-1}=du_m+cu_{m-1},\quad \beta_{m-1}=-bu_m-au_{m-1}.
$$
Arguing by induction, one sees that there exist
elements $u_{m-2},\ldots,u_0\in \mathcal{K}[X_1,\ldots,X_n]$ such that
$$
\alpha_{i}=du_{i+1}+cu_{i},\quad \beta_{i}=-bu_{i+1}-au_{i},\quad i=0,\ldots,m-1.$$
Now,
\begin{eqnarray*}
AF+BG=\alpha_0b+\beta_0d&=&(du_1+cu_0)b+(-bu_1-au_0)d\\
&=&(bc-ad)u_0=-\rho u_0.
\end{eqnarray*}
Thus, $AF+BG\in (\rho)$ and hence $\mathcal{J}\subset (\rho)$.
Since conversely $\rho=aG-cF\in\mathcal{J}$, we find $\mathcal{J}=(\rho)$. \CVD

We recall the definition (\ref{defpartial}) of the operators $\partial_n^{(d)}$:
\begin{equation}\label{defpartial2}
\partial_n^{(d)}f = D_nf + \sum_{i=1}^n (-1)^i{d+n-1\choose
  i}(D_{n-i}f)(D_{i-1}E).
\end{equation}
If $f\in\widetilde{M}^{\leq l}_{w,m}\setminus
  \widetilde{M}^{\leq l-1}_{w,m}$,
we write $\partial_nf$ instead of $\partial^{(w-l)}_nf$.

\begin{Lemma}\label{lemmepartialab}
Let $k\geq 1$. Let $f=aE+b$ be in
$\widetilde{M}^{\leq 1}_{q^k+1,1}\setminus M$,
with $a,b\in M$, and such that 
$f$ divides $D_1f$. Then $f=-D_1a$ and $D_1f=0$.
\end{Lemma}
\noindent\emph{Proof.}
We have $D_1f=Af$ with $A\in\widetilde{M}^{\leq 1}_{2,1}$, so that
$D_1f=\mu Ef$ with $\mu\in C$. By definition (\ref{defpartial2})
we have $D_1f=d(f)Ef+\partial_1f=\partial_1f$, hence
$\partial_1f=\mu Ef$. Since 
by Theorem~\ref{proposition_partial} $\partial_1f$
has depth $\leq 1$, we deduce $\mu=0$.
Definition (\ref{defpartial2}) gives 
$\partial_1a=D_1a+aE$ and $\partial_1b=D_1b-bE$.
Therefore,
\begin{eqnarray*}
0=D_1f=(D_1a)E+a(D_1E)+D_1b&=&(\partial_1a-aE)E+aE^2+\partial_1b+bE\\
&=&(\partial_1a)E+\partial_1b+bE.
\end{eqnarray*}
Since $\partial_1a,\partial_1b\in M$ by Theorem~\ref{proposition_partial},
we get
$\partial_1a=-b$ and $\partial_1b=0$. This yields
$f=aE-\partial_1a=-D_1a$. \CVD

\subsubsection{Proof of Theorem \ref{theorem:xkbis}.} 

The fact that $\epsilon_D(x_k)=e(k+1)$ follows from
Proposition \ref{propformulasgsxs}.

The theorem is clear for $k=0$ because $\widetilde{M}^{\leq 1}_{2,1}$
has dimension $1$ and is generated by $E=-x_0$.

For $k>0$, let $f=aE+b$, with $a,b\in M$, be normalised in $\widetilde{M}^{\leq 1}_{q^k+1,1}\setminus M$ and differentially extremal, so that 
$\epsilon_D(f)\geq e(k+1)$. We want to prove that $f$ is proportional to $x_k$.
Let us assume by contradiction that this is not the case.
We will use Lemma \ref{propocornelissen}, which 
forces us to distinguish the even and the odd $k$ cases.

\medskip

\noindent\emph{Case of $k$ odd.} We claim that $a,g_{k+1}$ are coprime
in $M$. Let us assume the contrary.
We first note that $f$, thus $a$, is
defined over the inseparable closure of $K$ by Lemma \ref{inseparable}.
Now, since $k+1$ is even, we can apply Lemma
\ref{propocornelissen} to $g_{k+1}$. We deduce that
if $a,g_{k+1}$ have a non-constant common factor, then
it must be $g_{k+1}$ (up to a non zero constant).
So $g_{k+1}$ must divide $a$ which is obviously false,
by elementary weight consideration.

The ideal ${\mathcal I}_k\cap M$ of $M=C[g,h]$ is principal,
generated by the resultant $\rho=\mathbf{Res}_E(f,x_{k+1})$
by Lemma \ref{lemmaresultant}. This is a modular
form of weight $w=q^{k}(q+1)$, by Lemma \ref{lemme_resultant},
and it is non-zero since $x_{k+1}$ and $f$ are coprime
by the irreducibility of $x_{k+1}$ (Lemma \ref{xkirreducible}).

We have $D_jf/f\in\widetilde{M}$ and $D_jx_{k+1}=0$
for $1\leq j\leq q^{k+1}-1$.
Hence, the ideal ${\mathcal I}_k=(f,x_{k+1})$ of $\widetilde{M}$
contains the sets $D_1{\mathcal I}_k,\ldots,D_{q^{{k+1}}-1}{\mathcal I}_k$.

By Theorem~\ref{proposition_partial},
$\partial_j\rho\in{\mathcal I}_k\cap M,\quad j=1,\ldots,q^{k+1}-1$. Since
${\mathcal I}_k\cap M$ is principal, we get:
\begin{equation}(\partial_j\rho)/\rho\in M,\quad j=1,\ldots,q^{k+1}-1.
\label{eq:partialnw}\end{equation}

It is then straightforward to see, from Identity (\ref{defpartial2}), that
\begin{equation}(D_1\rho)/\rho,\ldots,(D_{q^{k+1}-1}\rho)/\rho\in
\widetilde{M}.\label{Didividesrho}\end{equation}

By \cite[Lemma 5.7]{BP} and Proposition \ref{f0hs},
the property (\ref{Didividesrho}) implies that
$$\rho=h^{\nu_\infty(\rho)}\phi^{q^{k+1}},\quad \mbox{ for some $\phi\in M$ with }\nu_\infty(\phi)=0.$$ This equality of modular forms yields the following equality of weights:
\begin{equation}\label{eq:equalityofweights}
q^k(q+1)=\nu_\infty(\rho)(q+1)+w(\phi)q^{k+1}.\end{equation}

If $q>2$ (resp. $q=2$), the equality above can hold only if $w(\phi)=0$ (resp.
$w(\phi)=0,1$). Indeed, when $w(\phi)\neq 0$, we have $w(\phi)\geq q-1$.
But the non-negativity of $\nu_\infty(\rho)$ in (\ref{eq:equalityofweights}) is 
contradictory with the conditions $q>2$ and $w(\phi)\geq q-1$, or 
$q=2$ and $w(\phi)\geq 2$.

Let us assume, for $q\geq 2$, that $w(\phi)=0$, that is, that
$\phi$ is a non-zero constant. There exists $c$ in $C^{\times}$,
such that $$\rho=-g_{k+1}f+ax_{k+1}=ch^{q^k}.$$
We have proved in Proposition \ref{xkordqk} that
$\nu_\infty(x_{k+1})=q^{k+1}$.
It then follows that $\nu_\infty(f)=\nu_\infty(g_{k+1}f)=q^k$ (remember that $\nu_\infty(g_{k+1})=0$). By the extremality of $x_k$,
$f$ is proportional to $x_k$,
a contradiction. 

It remains to treat the case $q=2$ and $w(\phi)=1$ in (\ref{eq:equalityofweights}). This yields $3\cdot 2^{k}=3\nu_\infty(\rho)+2^{k+1}$,
that is, $3\nu_\infty(\rho)=2^k$, which is impossible by the integrality of $\nu_\infty(\rho)$. 
This concludes the proof of the theorem in the case $k$ odd.

\medskip

\noindent\emph{Case of $k$ even.} In this case,
with the same arguments that we have used for the case $k$ odd,
we see that two subcases hold. The first is the case $f=aE+b$ with 
$a,g_k$ coprime, the second is the case of $a$ proportional to $g_k$.

\medskip

\noindent\emph{First subcase.} This case can be handled in about the same way 
as the case of $k$ odd, again with additional nested subcases corresponding to $q>2$ and $q=2$. Regardless to the value of $q$, the resultant 
$\rho=\mathbf{Res}_E(f,x_k)$ has weight $2q^k$ and satisfies:
$$\rho=h^{\nu_\infty(\rho)}\phi^{q^{k+1}},$$
by means of arguments very similar to that we have used in the case $k$ odd.
This yields
\begin{equation}\label{eq:equalityofweightscaseeven}
2q^k=\nu_\infty(\rho)(q+1)+w(\phi)q^{k+1}.\end{equation}
From this we see that $2q^k\geq w(\phi)q^{k+1}$, so that
Identity (\ref{eq:equalityofweightscaseeven}) can hold only when
$w(\phi)=0$ if $q>2$ or when $w(\phi)=0,1$ if $q=2$.
But $w(\phi)=0$ is impossible in (\ref{eq:equalityofweightscaseeven})
since $q+1$ never divides $2q^k$. 
It remains to treat the case $q=2$ and $w(\phi)=1$.
In this case, $\nu_\infty(\rho)=0$, $\phi$ is proportional to $g$,
and there exists $c$ in $C^\times$
such that $\rho=-g_{k}f+ax_{k}=cg^{q^{k+1}}$.
By Proposition \ref{f0hs}, $D_1f=\cdots=D_{q^{k+1}-1}f=0$.
However, by \cite[Proposition 5.4, Identity (37)]{BP}, the depth of $D_{q^k}f$
is equal to $q^k+1$ and $D_{q^k}f$ cannot vanish; a contradiction.

\medskip

\noindent\emph{Second subcase.} 
It remains to treat the case of $f=aE+b$ with $a$ proportional to $g_k$. But in this case,
by Lemma \ref{lemmepartialab}, $f$ is proportional to $D_1g_k$, thus proportional
to $x_k$ by definition of the sequence $(x_k)_{k\geq 0}$. The proof of the theorem is now complete. \CVD

\subsection{Numerical observations\label{experiments}}

We made several numerical computations, essentially with $q=2,3,5$,
thanks to the algorithms of Section \ref{algorithm}
and with the help of a computer.
Let us introduce the following family of quasi-modular
forms (where we recall that $L_k:=[k]\cdots[1]$ if $k> 0$ and $L_0:=1$):
\begin{equation*}
\eta_k:=L_{k+1}x_k^q+L_k^qgx_{k+1}\in\widetilde{M}^{\leq q}_{q(q^k+1),1}
\setminus\widetilde{M}^{\leq q-1}_{q(q^k+1),1},\quad k\geq 0.
\end{equation*}
It is easy to prove that $\nu_\infty(\eta_k)=q^{k+1}+q-1$ and 
$\epsilon_D(\eta_k)=0$, for all $k\geq 0$ (use Lemma \ref{lemma:epsilonD} for the latter identity).

We present a table describing both the results
of this paper and the analysis of the results of the numerical
experiments we made. We look at extremal and differentially extremal
quasi-modular forms in the vector 
spaces $\widetilde{M}^{\leq l}_{l(q^k+1),l}$ with
$l=1,\ldots,q+1$ and $k\geq 0$. The integer $s$ is supposed to
be $<\frac{q+1}{2}$
while the integer $s'$ satisfies $\frac{q+1}{2}\leq s'\leq q-1$.

In the first column we enter the integer $l$ which determines the vector
space where we  look for an extremal quasi-modular
form. In the second column, there is an extremal quasi-modular form $f$
(unique up to multiplication by a scalar in $C^\times$).
If an asterisk $(*)$ appears, the resulting form is certified by
numerical computations solely in the case of $q=2,3,5$ and $k=0,1$.
If no asterisk figures, the validity of the entry is understood for all $q$
and for all $k\geq 0$. 

The third column contains $\nu_\infty(f)$.
In the fourth column, the differential exponent $\epsilon_D(f)$
of the corresponding form $f$ is computed or estimated. Again, in absence
of an asterisk, the result is unconditional. Otherwise,
the value of the entry represents a lower bound for $\epsilon_D(f)$
and its exact value for $q=2,3,5$ and $k=0,1$.

In the fifth column we enter a differentially extremal quasi-modular
form $f'$ in $\widetilde{M}^{\leq l}_{l(q^k+1),l}\setminus
\widetilde{M}^{\leq l-1}_{l(q^k+1),l}$. In this case, there is no
need to write down all the orders of vanishing 
and differential exponents, since they can be easily computed applying
the results of this text and in particular
Lemma \ref{lemma:epsilonD}. Of course, the presence of the asterisk
tells that the corresponding quantity is checked only for $q=2,3,5$
and $k=0,1$. It seems that these differentially extremal forms are
unique up to multiplication by a scalar in $C^\times$.

\bigskip

\begin{tabular}{|l|ll|ll|ll|ll|}
\hline
$l$ & $f$ & & $\nu_\infty(f)$ & & $\epsilon_D(f)$ & & $f'$ & \\
\hline
$1$ & $x_k$ & & $q^k$ & & $(k+1)e$ & & $x_k$ & \\
$2$ & $x_k^2$ & & $2q^k$ & & $(k+1)e$ & & $x_k^2$ & * \\
$\vdots$ & $\vdots$ & & $\vdots$ & & $\vdots$ & & $\vdots$ & \\
$s$ & $x_k^s$ & & $sq^k$ & & $(k+1)e$ & & $x_k^s$ & * \\
\hline
$s'$ & $x_k^{s'}$ &*  & $s'q^k$ & & $(k+1)e$ & & $x_k^{s'}$ & * \\
$\vdots$ & $\vdots$ & & $\vdots$ & & $\vdots$ & & $\vdots$ & \\
$q-1$ & $x_k^{q-1}$ &* & $q^{k+1}-q^k$ & & $(k+1)e$ & & $x_k^{q-1}$ & * \\
\hline
$q$ & $\eta_k$ &* & $q^k+q-1$ & & $0$ & & $x_k^q$ & * \\
\hline
$q+1$ &$\xi_k$ &* & $q^{k+2}+q^k$ & & $(k+2)e$ &*& $\xi_k$ & * \\
\hline
\end{tabular}

\bigskip

This table indicates that extremality and differential extremality are inequivalent conditions.
However, it seems that certain forms, notably the $x_k$'s and $\xi_k$'s, have the interesting
property of being at once irreducible, extremal and differentially extremal. 

These {\em primitive} forms will become, in the opinion of the authors, of particular 
importance and will deserve a crucial role in the forthcoming researches in this topic.
In the case $q=2$, computations have been pushed forward to higher values of $l$ and $k$
disclosing the existence of other primitive forms, that will be studied elsewhere.

\subsection{The forms $x_k$'s as solutions of certain differential systems.\label{hypergeometric}}


In this section, we give yet another property that characterises the collection of forms $x_k$ up to 
scalars in $C^\times$. The theorem below could be understood as a reasonable substitute 
of Theorem 2 of \cite{KK2}. We recall that the operators $\partial^{(d)}_n$
have been defined in (\ref{defpartial}).

\begin{Theorem}\label{corollary:hypergeometric}
Let $f$ be a non-constant element of $\widetilde{M}^{\leq 1}_{q^k+1,1}$. We have $$\partial^{(q^k)}_1f=\cdots=\partial^{(q^k)}_{q^{k+1}-1}f=0$$
if and only if $f$ is proportional to $x_k$.
\end{Theorem}

The proof of this theorem requires the two following lemmas:

\begin{Lemma}
We have:
\begin{equation}\label{Dqk1E}
D_{iq^k-1}E=E^{iq^k},\quad k\geq 0,\quad  1\leq i <q.
\end{equation}
\end{Lemma}
\noindent\emph{Proof.} It follows by the same arguments of the proof
of \cite[Lemma 3.8]{BP}. Since we have \cite[p. 686]{Ge}
$$ E = \frac{1}{\widetilde{\pi}}\sum_{a\in A\text{ monic }}\sum_{b\in A}
\frac{a}{az+b},$$
we deduce
\begin{eqnarray*}
D_{iq^k-1}E&=&\frac{1}{\widetilde{\pi}^{iq^k}}\sum_{a\in A\text{ monic }}
\sum_{b\in A}\left(\frac{a^i}{(az+b)^i}\right)^{q^k}\\
&=&(D_{i-1}E)^{q^k}\\
&=&E^{iq^k}.
\end{eqnarray*}\CVD

\begin{Lemma}\label{lemme:partial=0xk}
Let $f$ be a non-constant element of $\Psi_{e(k+1)-1}$ and $n$ be an integer $\leq q^{k+1}-1$. The following conditions are
equivalent:
\begin{enumerate}
\item $D_nf=0$ if $q^k$ does not divide $n$ and $D_nf=E^{iq^k}f$ if 
$n=iq^k$ with $0\leq i<q$.
\item $\partial^{(q^k)}_1f=\cdots=\partial^{(q^k)}_{q^{k+1}-1}f=0$.
\end{enumerate}
\end{Lemma}
We first need to focus on certain structures which appear in the operator $\partial_{n}^{(q^k)}$.
With $d=q^k$ and $1\leq n\leq q^{k+1}-1$, the definition of the operator $\partial^{(d)}_n$
reads:
$$\partial_n^{(q^k)}f=D_nf+R_1+R_2,$$
where 
\begin{eqnarray*}
R_1&=&(-1)^n\sum_{0\leq i\leq n-1,q^k\nmid i}(-1)^i\binomial{q^k+n-1}{n-i}(D_if)(D_{n-i-1}E),\\
R_2&=&(-1)^n\sum_{0\leq i\leq n-1,q^k|i}(-1)^i\binomial{q^k+n-1}{n-i}(D_if)(D_{n-i-1}E).
\end{eqnarray*}
We have the congruences modulo $p$:
\begin{equation}\label{congruencebinom1}
\binomial{q^k+n-1}{n-jq^k}\equiv\left\{\begin{array}{cc}\binom{i}{i-j} & \text{ if }n=iq^k\\
0 & \text{ if }q^k\nmid n\end{array}\right.,\quad \text{ for }\begin{array}{l}0\leq j\leq q-1,\\ jq^k\leq n\leq q^{k+1}-1.\end{array}
\end{equation}
The binomial in the left-hand side of the congruence is, up to multiplication by a power of $-1$,
the coefficient of $(D_{jq^k}f)(D_{n-jq^k-1}E)$ in the expression defining $\partial^{(q^k)}_nf$.
Hence, after (\ref{congruencebinom1}) and (\ref{Dqk1E}), 
\begin{eqnarray}
R_2&=&(-1)^i\sum_{j=0}^{i-1}(-1)^j\binom{i}{i-j}(D_{jq^k}f)(D_{(i-j)q^k-1}E)\nonumber\\
&=&(-1)^i\sum_{j=0}^{i-1}(-1)^j\binom{i}{i-j}(D_{jq^k}f)E^{(i-j)q^k},\quad \text{ if }n=iq^k,\label{case1}
\end{eqnarray}
and
\begin{equation}
R_2=0\quad\text{ otherwise. }\label{case2}
\end{equation}
\noindent\emph{Proof of Lemma \ref{lemme:partial=0xk}. 1 $\Rightarrow$ 2.} 
We prove that $\partial_n^{(q^k)}f=0$ for $n=1,\ldots,q^{k+1}-1$.
We have $R_1=0$ by hypothesis. If $q^k\nmid n$, $D_nf=0$ and $\partial_n^{(q^k)}f=0$ because of (\ref{case2}).
If $n=iq^k$, $D_nf=E^{iq^k}f$ by hypothesis and by (\ref{case1})
\begin{eqnarray}
R_2&=&(-1)^iE^{iq^k}f\sum_{j=0}^{i-1}(-1)^j\binomial{i}{i-j}\\
&=&-E^{iq^k}f.\label{eqR2}
\end{eqnarray}
Hence, in this case too, $\partial_n^{(q^k)}f=D_nf+R_2=0$.

\medskip

\noindent\emph{2 $\Rightarrow$ 1.} For $n=1$, the statement is true because $\partial^{(q^k)}_1f=D_1f$.
Let us assume that 
we have already proved that, for all $n\leq m-1$ ($m$ being an integer $\geq 2$)
$D_nf=0$ if $q^k$ does not divide $n$ and $D_nf=E^{iq^k}f$ if 
$n=iq^k$. We know that $\partial_m^{(q^k)}f=0$. Hence
$D_mf=-R_1-R_2=-R_2$ (because $R_1=0$ by induction hypothesis).
If $q^k\nmid m$, $R_2=0$ by (\ref{case2}) and $D_mf=0$.
Otherwise, the induction hypothesis implies Equality (\ref{eqR2}) and 
$D_mf=E^mf$.\CVD

\noindent\emph{Proof of Theorem \ref{corollary:hypergeometric}.} If $f=cx_k$, $c\in C^\times$, then Identity (\ref{xixs})
of Proposition \ref{propformulasgsxs} 
and Lemma \ref{lemme:partial=0xk}
imply that the identities involving the operators $\partial^{(q^k)}_n$ hold.
For the other implication,  we first apply
Lemma~\ref{lemme:partial=0xk} and then, Theorem \ref{theorem:xkbis}.\CVD

\begin{Remark}{\em We could not find an analogue of Theorem \ref{corollary:hypergeometric}
for the family $(\xi_k)_{k\geq 0}$. Similarly, we did not find a reasonable substitute of Theorem 1 of \cite{KK2}.}\end{Remark}

\subsection{Structure of the subfields $F_k$'s\label{Fk}}

The content of this section is independent on our quest of finding differentially extremal quasi-modular forms
and can be skipped in a first reading of the paper.

Proposition \ref{f0hs} says that every differentially extremal quasi-modular form is multiple by a power of the form $h$
of a (necessarily isobaric) element of $F_k$. By \cite[Proposition 2.2]{MP}, there exists, for all $k\geq 0$, an element
$z_k\in F\setminus F_k$ such that $F$ is a $F_k$-vector space of dimension $p^k$ of basis $(1,z_k,\ldots,z_k^{p^k-1})$ and such that
for all $n,m\in\NN$, $D_nz_k^m=\binom{m}{n}z_k^{m-n}$. However, this does not clarify much the structure of 
the fields $F_k$ themselves.

We now define the sequence
$(y_k)_{k\ge 1}$ by $$y_k=\Delta^{q^{k-1}}x_{k-1}.$$

Here, we prove: 

\begin{Theorem}\label{generators}
Let $r\ge 0$ be a non-negative integer. For every integer $k$ such that
$q^{k+1}>p^r$ we have
$$ F_r=F_{r-1}^p(x_{k+1},y_{k+1}).$$
\end{Theorem}

For instance, if we take $k=r$ we find $ F_k=F_{k-1}^p(x_{k+1},y_{k+1})$
for all $k\ge 0$.


To prove Theorem~\ref{generators} we will need two lemmas.

\begin{Lemma}\label{index}
For all $r\ge 0$ we have
$$ [F_r:F_{r-1}^p]=p^2. $$
\end{Lemma}

\noindent\emph{Proof.} We know that $[F_{r-1}:F_r]=p$ for all $r$.
Let us prove by induction that
$[F_r:F_r^p]=p^3$. This is clear for $r=-1$. If now this property
holds for the integer $r-1$, then:
$$ [F_r:F_r^p]=\frac{[F_{r-1}:F_{r-1}^p][F_{r-1}^p:F_r^p]}{[F_{r-1}:F_r]}=
\frac{p^3\cdot p}{p}=p^3,$$
hence the property also holds for $r$. Thus it holds for all $r\ge 0$,
from which we deduce
$$ [F_r:F_{r-1}^p]=\frac{[F_r:F_r^p]}{[F_{r-1}^p:F_r^p]}=\frac{p^3}{p}=p^2. $$
\CVD

\begin{Lemma}\label{independence}
For all $k\ge 0$ the elements $x_{k+1}^iy_{k+1}^j$ ($1\le i,j\le p-1$)
are linearly independent over $F^p$.
\end{Lemma}

\noindent\emph{Proof.} Let us first consider the case $k=0$.
Suppose that there exists a relation
\begin{equation}\label{aze}
\sum_{\genfrac{}{}{0pt}{2}{0\leq i\leq p-1}{0\leq j\leq p-1}}
\lambda_{ij}^px_1^iy_1^j=0
\end{equation}
where $\lambda_{ij}\in F$.
Since $x_1=-(Eg+h)$ and $y_1=Eh^{q-1}$ this rewrites, after multiplying
(\ref{aze}) by $h^{p-1}$:
\begin{equation}\label{azer}
\sum_{\genfrac{}{}{0pt}{2}{0\leq i\leq p-1}{0\leq j\leq p-1}}
\bigl((-1)^{i}\lambda_{ij}h^{jp^{e-1}}\bigr)^p(Eg+h)^iE^jh^{p-1-j}=0.
\end{equation}
Since the forms $Eg+h$, $E$ and $h$ are algebraically independent over
$C$ by \cite[Lemma~2.4]{BP}, the forms $(Eg+h)^iE^jh^k$ ($0\le i,j,k\le p-1$)
are obviously linearly independent over $F^p$. In particular,
it follows from (\ref{azer}) that $(-1)^i\lambda_{ij}h^{jp^{e-1}}=0$
for all $i,j$, hence $\lambda_{ij}=0$ for all $i,j$.
This proves the result for $k=0$.

Suppose now that $k\ge 1$. Since $x_{k+1}=g^{q^k}x_k-
[k]y_k$ and $y_{k+1}=\Delta^{q^{k}}x_{k}$, we have:
$$ F^p(x_{k+1},y_{k+1})=F^p(x_{k+1},x_k)=
F^p(g^{q^k}x_k-[k]y_k,x_k)=F^p(x_k,y_k).$$
By induction, it follows that $F^p(x_{k+1},y_{k+1})=F^p(x_1,y_1)$,
hence
$$[F^p(x_{k+1},y_{k+1}):F^p]=[F^p(x_1,y_1):F^p]=p^2$$
by the case  $k=0$ we have just proved. But this means that
the $p^2$ generators $x_{k+1}^iy_{k+1}^j$ ($1\le i,j\le p-1$) of the
$F^p$-algebra $F^p(x_{k+1},y_{k+1})$ are $F^p$-linearly independent.
\CVD

\noindent\emph{Proof of Theorem~\ref{generators}.}
Let $r\ge 0$ and $k\ge 0$ be as in the theorem.
First of all, we note that by Lemma~\ref{lemme:2ek} we have
$x_{k+1}\in F_r$. Next, the relation $[k+1]y_{k+1}=x_{k+1}g^{q^{k+1}}-x_{k+2}$
together with the fact that $g^{q^{k+1}}$, $x_{k+1}$, $x_{k+2}$ all
belong to $F_r$ shows that $y_{k+1}\in F_r$. Thus we have the inclusions
$$ F_{r-1}^p\subset F_{r-1}^p(x_{k+1},y_{k+1})\subset F_r.$$
Now, we have $[F_r:F_{r-1}^p]=[F_{r-1}^p(x_{k+1},y_{k+1}):F_{r-1}^p]=p^2$
by Lemmas \ref{index} and \ref{independence}. It follows that
$F_r=F_{r-1}^p(x_{k+1},y_{k+1})$.\CVD

\begin{Remark}{\em 
Let us define, for $s\in\ZZ$ and $k\geq 0$: $$\mathcal{A}_k^{(s)}:=\{f\in\Psi_k\text{ such that }\nu_\infty(f)\equiv s\pmod{p^{k+1}}\}\cup\{0\}.$$
It is an $F_k$-vector space by Proposition \ref{f0hs} and \cite[Lemma 5.2]{BP}.
Then it is easy to see that we have the following direct sum:
$$\mathcal{A}_k:=C[\Psi_k]=\bigoplus_{s\in\ZZ/p^{k+1}\ZZ}\mathcal{A}_k^{(s)}.$$
$\mathcal{A}_k$ is a $\ZZ/p^{k+1}\ZZ$-graded $F_k$-algebra.
Proposition \ref{f0hs} implies that a basis of this algebra is $(1,h,\ldots,h^{p^{k+1}-1})$.
The difficulty of constructing differentially extremal quasi-modular
forms comes from the 
difficulty of computing the intersections $\mathcal{A}_k^{(s)}\cap\widetilde{M}^{\leq l}_{w,m}$ for given $s,k,l,w,m$.
This seems to explain why we did not really take advantage of
Theorem \ref{generators}.

It is easy to deduce from Proposition~\ref{f0hs} that
$\mathcal{A}_k^{(s)}\cap M=M^{p^{k+1}}h^s$ with $M=C[g,h]$, that is,
Lemma 5.7 of \cite{BP} (this result also follows from  \cite[Theorem 2.6]{ZR}).
}\end{Remark}

\section{Differential tools\label{differential_tools}}

This section, divided in two distinct subsections, contains two contributions to the study of differential properties of Drinfeld quasi-modular forms,
the proof of Theorem \ref{proposition_partial} for the operators $\partial^{(d)}_n$ and 
the description of an algorithm which allows to compute higher derivatives of Drinfeld quasi-modular forms.
Although these tools have been already used, namely in Section \ref{differential_extremality}, we decided to
collect them in a separate section as they can be of interest independent on the study of differential extemality. 

\subsection{Higher Serre's operators.\label{sserre_operators}}

Here we study the higher Serre's operators (\ref{defpartial}) and prove
Theorem~\ref{proposition_partial}.

Let $n,d$ be non-negative integers. We have defined the {\em $n$-th Serre's
operator of degree $d$}, $\partial_n^{(d)}:\widetilde{M}\rightarrow
\widetilde{M}$, by the following formula, that we quote again to ease the reading of this section:
\begin{equation}\label{defpartial3}
\partial_n^{(d)}f = D_nf + \sum_{i=1}^n (-1)^i{d+n-1\choose i}(D_{n-i}f)(D_{i-1}E).
\end{equation}

Notice that in this definition the integer $d$ is arbitrary:
In particular, it is not necessarily the weight or the degree of $f$
(at this stage $f$ is not supposed to be a quasi-modular form, anyway).
If $n=0$, then $\partial_n^{(d)}f=\partial_0^{(d)}f=f$.
If $n=1$, we have
$\partial_n^{(d)}f=D_1f-dEf$, which coincides with the formula
(\ref{defserreop}) when $d=w(f)$.

In \cite{KK1}, the authors use the family of operators $\theta_k^{(r)}$, described in the introduction, which act on vector spaces
of quasi-modular forms for $\mathbf{SL}_2(\ZZ)$. When $n\geq 1$, there is a strong similarity between 
our operator $\partial_n^{(d)}$ and the operator $\frac{1}{n!}\theta_d^{(n-1)}$.

In \cite{US}, the authors introduce a class of operators on Drinfeld modular forms which 
could play the role of Rankin-Cohen operators in the drinfeldian framework. 
Their definition appears in formula (3.14) of their Theorem 3.7, and the notation they adopt for their operator is $[\cdot,\cdot]_{k,l,n}$.
It is easy to prove that for all $n,d\geq 1$ there exists $\lambda_{n,d}\in\FF_q$ such that:
$$\partial_n^{(d)}f=D_nf-\lambda_{n,d}(-\widetilde{\pi})^{1-n}(n+d-1)[E,f]_{2,d,n-1}.$$
However, we do not have a general receipt to compute $\lambda_{n,d}$; for instance, we notice that it can vanish
for certain choices of $n,d$.

The remarkable feature of the operator $\partial_n^{(d)}$ is that it
does not increase the depth of quasi-modular forms of degree $d$. In the
classical case, a proof of this fact appears, for example, in
\cite[Proposition 3.3]{KK1}.

\begin{Theorem}\label{proposition_partial}
Let $w,l$ be non negative integers with $w\geq 2l$ and let $m$ be a
class in $\ZZ/(q-1)\ZZ$. Define $d:=w-l$. For every $n\in\NN$, we have
\begin{equation}\label{eq:partial}
\partial_n^{(d)}(\widetilde{M}^{\leq l}_{w,m})\subset \widetilde{M}^{\leq l}_{w+2n,m+n}.
\end{equation}
In particular, $\partial_n^{(w)}$ sends modular forms of weight $w$
and type $m$ on modular forms of weight $w+2n$ and type $m+n$.
\end{Theorem}

To prove this Theorem we will use
the notion of polynomial associated with a Drinfeld quasi-modular form
introduced
in \cite{BP}. We recall for convenience the
definitions and properties we will need here. If $f$ is a
Drinfeld quasi-modular form of weight $w$ and type
$m$, then there exists a unique polynomial
$P_f=\sum_{i=0}^lf_iX^i\in\widetilde{M}[X]$
such that
$$
f(\gamma(z))= \frac{(cz+d)^w}{(\det\gamma)^m}\sum_{i=0}^{l}
f_i(z)(\frac{c}{cz+d})^i
$$
for every $z\in\Omega$ and every $\gamma=
\left(\begin{matrix}
a&b\cr
c& d\cr
\end{matrix}\right)
\in\mathbf{GL}_2(A)$.
This polynomial is the {\em polynomial associated with $f$}, its
coefficients are in fact Drinfeld quasi-modular forms. If $f=0$, then
$P_f=0$ and if $f\not=0$, then the
degree of $P_f$ is equal to the depth of $f$. When $f$ is a modular
form, then we have $P_f=f$ and for $f=E$ we have
$P_E=E-\widetilde{\pi}^{-1}X$. Finally, if  $f_1$, $f_2$ are two
Drinfeld quasi-modular forms, then $P_{f_1f_2}=P_{f_1}P_{f_2}$ and, if
$f_1,f_2$ have the same weight and the same type, then $P_{f_1+f_2}=P_{f_1}+P_{f_2}$.

It will be convenient to introduce the following notation, where 
$X$ is an indeterminate over the ring $\widetilde{M}$, which is equal to $C[E,g,h]$
in virtue of \cite[Theorem 1]{BP}.

\begin{Definition}\emph{
Let $f=f(E,g,h)$ be an element of $\widetilde{M}$. For every $n\in\NN$
we define $\partialbis_nf\in\widetilde{M}$ by the formula
\begin{equation}\label{defdelta}
f(E+X,g,h) = \sum_{n\geq 0} (\partialbis_nf)X^n\in C[E,g,h,X],
\end{equation}}
\end{Definition}

The family $\partialbis=(\partialbis_n)_{n\geq 0}$ is obviously an iterative higher
derivation on $\widetilde{M}$. We have
$\partialbis_1=\frac{\partial}{\partial E}$ and
$\bigl(\frac{\partial}{\partial E}\bigr)^n=n!\, \partialbis_n$.
The interest of introducing $\partialbis$ is that the following
property clearly holds (for all Drinfeld quasi-modular forms $f\not=0$
and all integers
$l\geq 0$):
\begin{equation}\label{equivalence}
l(f)\leq l \quad \Longleftrightarrow \quad \partialbis_nf=0 \ {\rm for\ all}\ n\geq l+1.
\end{equation}

In fact, it turns out that the $\partialbis_nf$'s are, up to powers of $-\widetilde{\pi}$,
the coefficients of the associated polynomial $P_f$.

\begin{Lemma}\label{lemme:delta_and_poly}
Let $f$ be a Drinfeld quasi-modular form. Then we have
$$
P_f(X)=\sum_{n\geq 0}\frac{\partialbis_nf}{(-\widetilde{\pi})^n}X^n.
$$
\end{Lemma}

\noindent\emph{Proof.}
We may assume that $f$ is non zero. Let $w$, $m$ and $l$ denote
respectively the weight, the type and the depth of $f$. Let us write
$f = \sum_{i=0}^lf_iE^i$,
where $f_i=f_i(g,h)\in M_{w-2i,m-i}$. We have
$$
P_f=\sum_{i=0}^{l}P_{f_iE^i}=\sum_{i=0}^lf_iP_E^i=\sum_{i=0}^lf_i(E-\frac{1}{\widetilde{\pi}}X)^i
=f(E-\frac{1}{\widetilde{\pi}}X,g,h),
$$
so the result immediately follows from (\ref{defdelta}).\CVD

This very simple result allows to translate any property
about the coefficients of the polynomial $P_f$ into a property about
the higher derivatives $\partialbis_nf$ and conversely. For example,
the reader will easily check that the fact that $\partialbis$ is
iterative is equivalent to Lemma~2.5 of \cite{BP}. For the proof of
Theorem~\ref{proposition_partial}, we will use the reformulation
of \cite[Proposition~3.1]{BP} in terms of the higher derivation $\partialbis$. This yields
the following commutation rule between $\partialbis_j$ and $D_n$.
Note that in the complex case and for $j=1$, the analogous formula
is established in \cite{KK1} (during the proof of Proposition~3.3).

\begin{Lemma}\label{lemme:commutation}
Let $j\geq 0$ and $n\geq 0$ be non negative integers, and let $f$ be a Drinfeld
quasi-modular form of weight $w$. Then we have
$$
\partialbis_jD_nf=\sum_{r=0}^n {w+n-j+r-1\choose r}D_{n-r}\partialbis_{j-r}f,
$$
where we set $\partialbis_i=0$ if $i<0$.
\end{Lemma}

\noindent\emph{Proof.}
This is exactly the formula for $P_{{\cal D}_nf}$ given in
Proposition~3.1 of \cite{BP}, expressed in terms of the higher
derivation $\partialbis$ with help of Lemma~\ref{lemme:delta_and_poly}, and
taking into account the fact that $D_i=\frac{1}{(-\widetilde{\pi})^i}{\cal D}_i$. \CVD

The next Lemma is analogous to a similar formula appearing in the
proof of \cite[Proposition~3.3]{KK1}. However, the proof in our case
cannot be done by induction as in \cite{KK1} and thus requires more care.

\begin{Lemma}\label{lemme:deltapartial}
Let $n$, $d$, $w$ and $k$ be non negative integers, and let $f$ be a Drinfeld
quasi-modular form of weight $w$. Then we have
\begin{equation}\label{eq:proppartial}
\partialbis_k\partial_{n}^{(d)}f=\sum_{i=0}^k(-1)^i{d+k-w-1\choose i}\partial_{n-i}^{(d+i)}
\partialbis_{k-i}f,
\end{equation}
where we define $\partial_j^{(d)}=0$ if $j<0$.
\end{Lemma}

\noindent\emph{Proof.}
The formula (\ref{eq:proppartial}) being plainly true if $n=0$, we will assume in the
following that $n\geq 1$.
Applying the operator $\partialbis_k$ to equation (\ref{defpartial3}) and using
Leibniz rule for the product $\partialbis_k((D_{n-i}f)(D_{i-1}E))$, we get
\begin{equation}\label{equation1}
\partialbis_k\partial_n^{(d)}f = \partialbis_kD_nf +
\sum_{\genfrac{}{}{0pt}{2}{1\leq i\leq n}{0\leq j\leq k}}
(-1)^i {d+n-1\choose i}(\partialbis_jD_{n-i}f)(\partialbis_{k-j}D_{i-1}E).
\end{equation}
Applying Lemma~\ref{lemme:commutation} with $f=E$ and noting that
$\partialbis_{k-j-r}E=0$ if $r>k-j$ or $r< k-j-1$, we have (we use the convention
$D_i=0$ if $i<0$)
\begin{eqnarray*}
\partialbis_{k-j}D_{i-1}E & = & \sum_{r\geq 0}{i-k+j+r\choose r}D_{i-1-r}\partialbis_{k-j-r}E\\
                  & = & {i \choose k-j}D_{i+j-k-1}E+{i-1\choose k-j-1}D_{i+j-k}(1)\\
                  & = & \left\{\begin{array}{cl}
                        {i \choose k-j} D_{i+j-k-1}E & {\rm if}\ j>k-i\cr
                        1 & {\rm if}\ j=k-i\cr
                        0 & {\rm if}\ j<k-i.\cr
                        \end{array}\right.
\end{eqnarray*}
Substituting in (\ref{equation1}), we obtain
\begin{equation}\label{aplusb}
\partialbis_k\partial_n^{(d)}f=A+B,
\end{equation}
where
\begin{equation*}\label{defA}
A = \sum_{0\leq i\leq n}(-1)^i{d+n-1\choose i}\partialbis_{k-i}D_{n-i}f
\end{equation*}
and
$$ B= \sum_{\genfrac{}{}{0pt}{2}{1\leq i\leq n}{k-i+1\leq j\leq k}}
(-1)^i{d+n-1\choose i}{i\choose k-j}(\partialbis_jD_{n-i}f)(D_{i+j-k-1}E).$$
Applying Lemma~\ref{lemme:commutation} again
and then making the change of variable $I=i+r$, we first find the following expression
for $A$:
\begin{eqnarray}
A &=& \sum_{0\leq i\leq n}\sum_{0\leq r\leq n-i}(-1)^i{d+n-1\choose i}
{w+n-k+r-1\choose r}D_{n-i-r}\partialbis_{k-i-r}f\nonumber\\
 &=& \sum_{0\leq I\leq n}\left[\sum_{0\leq i\leq I}(-1)^i{d+n-1\choose i}
{w+n-k+I-i-1\choose I-i}\right] D_{n-I}\partialbis_{k-I}f\nonumber\\
&=& \sum_{0\leq I\leq n} (-1)^I{d+k-w-1\choose I}D_{n-I}\partialbis_{k-I}f.\label{aaa}
\end{eqnarray}
Similarly for $B$, we apply Lemma~\ref{lemme:commutation} to $\partialbis_jD_{n-i}f$
and then make the changes of variable $J=i+j-k$, $I=i+r-J$. This yields
\begin{equation*}
B= \sum_{\genfrac{}{}{0pt}{2}{0\leq I\leq n-1}{1\leq J\leq n-I}}(-1)^{I+J}S_{IJ}
(D_{n-I-J}\partialbis_{k-I}f)(D_{J-1}E),
\end{equation*}
where
\begin{equation*}
S_{IJ}=\sum_{r=0}^I(-1)^r{d+n-1\choose I+J-r}{I+J-r\choose I-r}{w+n-k-J+r-1\choose r}.
\end{equation*}
Since
$${d+n-1\choose I+J-r}{I+J-r\choose I-r}={d+n-1\choose d+n-1-J}{d+n-J-1\choose I-r},$$
the sum $S_{IJ}$ is equal to
$$ {d+n-1\choose d+n-1-J} T_{IJ}$$
with
\begin{eqnarray*}
T_{IJ} & = & \sum_{r=0}^I(-1)^r{d+n-J-1\choose I-r}{w+n-k-J+r-1\choose r}\\
& = & {d+k-w-1\choose I}
\end{eqnarray*}
by \cite[Lemma 3.2]{BP}.
Hence we obtain the following expression for $B$:
\begin{eqnarray}
B&=&\sum_{\genfrac{}{}{0pt}{2}{0\leq I\leq n-1}{1\leq J\leq n-I}}(-1)^{I+J} {d+k-w-1\choose I}
{d+n-1\choose d+n-1-J}\qquad\qquad \nonumber\\
&&\phantom{\sum_{\genfrac{}{}{0pt}{2}{0\leq I\leq n-1}{1\leq J\leq n-I}}(-1)^{I+J}(d+k-w-1)}
\times (D_{n-I-J}\partialbis_{k-I}f)(D_{J-1}E).\label{bbb}
\end{eqnarray}
Now, we note that in the formula (\ref{aaa}) the summands vanish if $I>n$ or $I>k$, so
we may assume that $I$ runs from $0$ to $k$. Similarly, in the formula (\ref{bbb})
we can let $I$ vary from $0$ to $k$. Using this remark, and substituting (\ref{aaa}) and
(\ref{bbb}) in (\ref{aplusb}), we get:
\begin{eqnarray*}
\partialbis_k\partial_n^{(d)}f & =
&\sum_{I=0}^k(-1)^I{d+k-w-1\choose I}\biggl[D_{n-I}\partialbis_{k-I}f\\
&& \phantom{sum} + \ \sum_{J=1}^{n-I}(-1)^J
{d+n-1\choose d+n-1-J}(D_{n-I-J}\partialbis_{k-I}f)(D_{J-1}E)\biggr].
\end{eqnarray*}
This is nothing else than the formula (\ref{eq:proppartial}), so Lemma~\ref{lemme:deltapartial}
is proved.\CVD

We can now prove Theorem~\ref{proposition_partial}.

\noindent\emph{Proof of Theorem~\ref{proposition_partial}.}
Let $f\in\widetilde{M}_{w,m}^{\leq l}$, and let $k$ be any integer such that
$k\geq l+1$. We want to show that $\partialbis_k\partial_n^{(d)}f=0$, which will
prove the proposition by Property (\ref{equivalence}).
We use for this Lemma~\ref{lemme:deltapartial}. Since
here $d=w-l$, we have
\begin{equation}\label{ccc}
\partialbis_k\partial_{n}^{(d)}f=\sum_{i=0}^k(-1)^i{k-l-1\choose i}\partial_{n-i}^{(d+i)}
\partialbis_{k-i}f.
\end{equation}
In this sum, if $i\leq k-l-1$, then $k-i\geq l+1$ and hence $\partialbis_{k-i}f=0$
(Property (\ref{equivalence})).
If now $i\geq k-l$, then ${k-l-1\choose i}=0$. So all the summands
in the right-hand side of (\ref{ccc}) vanish and
$\partialbis_k\partial_n^{(d)}f=0$.
\CVD

\subsubsection{Digression: an application to eigenforms of Hecke operators.\label{heckeoperators}}

Theorem~\ref{proposition_partial} has already been used in the proof
of Theorem~\ref{theorem:xkbis} (in Section~\ref{differential_extremality}). Another interesting
application of this theorem is that it can be used to construct {\em a priori} new eigenforms
for Hecke operators, from given ones.

Let ${\gp}=(P)$ be a non zero prime ideal of $A$, where $P$ is a monic
polynomial. Following \cite[\S~1.8]{Go2} or \cite[\S~7]{Ge}, we define, for any quasi-modular
form $f\in\widetilde{M}_{w,m}^{\leq l}$, $T_{\gp}f$ by the formula
\begin{equation}\label{hecke}
(T_{\gp}f)(z)=P^wf(Pz)+\sum_{b\in A\atop \deg_\theta b<\deg_\theta P}f(\frac{z+b}{P}),
\end{equation}
where we remark the dependence of this operator on the weight $w$. We also notice that 
there is no reason for $T_{\gp}f$ to lie in $\widetilde{M}$, except when we already know that $f\in M$.

\begin{Lemma}\label{lemme:TP}
If $f\in\widetilde{M}^{\leq l}_{w,m}$ is a quasi-modular form which is an {\em eigenform} for $T_{\gp}$
with eigenvalue $\lambda\in C$ (that is, such that $T_{\gp}f=\lambda f$), then for all $n\geq 1$,
 $D_nf\in\widetilde{M}^{\leq l+n}_{w+2n,m+n}$ also is an eigenform with eigenvalue
$\lambda P^n$.
\end{Lemma}
\noindent\emph{Proof.}
The function $\varphi:z\mapsto f(Pz)$ satisfies $(D_n\varphi)(z)=P^n(D_nf)(Pz)$
for all $n$, and the functions $f_b:z\mapsto f((z+b)/P)$ satisfy $(D_nf_b)(z)=P^{-n}(D_nf)((z+b)/P)$.
Since the weight of $D_nf$ is $w+2n$, we then see that
$$ T_{\gp}(D_nf)=P^nD_n(T_{\gp}f) \qquad \text{for all\ } n\ge 0,$$ where this time, the operator
$T_{\gp}$ in the left-hand side is defined as in (\ref{hecke}) but with $w$ replaced by $w+2n$.
It immediately follows that if $T_{\gp}f=\lambda f$ for some
$\lambda\in C$, then  $T_{\gp}(D_nf)=\lambda P^n(D_nf)$.\CVD

\noindent\emph{Example.} The normalised Eisenstein series $g_k$ defined in (\ref{def_gk})
(for $k\geq 1$) are modular eigenforms of weight $q^k-1$ and type $0$ of all the operators
$T_{\gp}$ with corresponding eigenvalue $P^{q^k-1}$, by
\cite[Proposition~7.2]{Ge}. Therefore, Lemma \ref{lemme:TP} says that $x_k=D_1g_k$
is eigenform of all the operators $T_{\gp}$ with corresponding eigenvalue
$P^{q^k}$ (if $k\geq 1$).

\medskip

\noindent\emph{Question.} Are the forms $\xi_k$ of Section \ref{sectionxik} eigenforms for all the Hecke operators?

\medskip

More interesting is the particular case in which $f$ is a modular form which is known to be 
an eigenform for all the Hecke operators, and $n$ is an integer such that $D_nf$ is 
again a modular form. The next lemma 
implies that when $d\geq 2$, there are infinitely many $n$'s for which $D_n=\partial_n^{(d)}$,
and thus every eigenform $f$ of $T_{\gp}$ of degree $\geq 2$ potentially yields
other eigenforms (these forms may be identically zero, but in many examples
they are not).

\begin{Lemma}\label{Lemma:dkn}
Let $d,k$ be non-negative integers with $d\geq 2$ and $p^k(p-1)\ge d-1$, and write
$n=1-d+p^{k+1}$. Then $\partial_n^{(d)}=D_n$.
\end{Lemma}
\noindent\emph{Proof.} By using \eg\ the formula (14) of \cite{BP},
one easily checks that $$ \binom{d+n-1}{i}\equiv 0 \pmod{p}\qquad \text{for all\ } i=1,\ldots,n.$$
The conclusion follows by applying these congruences to the formula defining 
the operators $\partial_n^{(d)}$.\CVD

If $f$ is a modular form of weight $d\geq 2$ and if $n$ is the integer of Lemma \ref{Lemma:dkn},
$D_nf$ is a modular form by Theorem~\ref{proposition_partial}. Note that for a given $f$, there
might be other choices of $n$ for which $\partial_n^{(d)}f=D_nf$, but it is not difficult to show
that there is at most one more choice for $n$ than the one in the Lemma, for which
$ \binom{d+n-1}{i}\equiv 0 \pmod{p}$ for all $i=1,\ldots,n.$

Example: We know \cite[Corollary (7.6)]{Ge} that $h$ is an eigenform of all the Hecke
operators. Here $d=q+1$, and we can take $n=p^{k+1}-q$ for every $k$
such that $p^k(p-1)\geq q$. Then we get an infinite family of
eigenforms $(D_{p^{k+1}-q}h)$.
We remark that we already knew from \cite[Lemma~3.10]{BP} that these
functions are modular. We do not know yet how to characterise the integers $n$ such that
$D_nh\neq 0$. 


By using the tools developed in Section \ref{algorithm}, some of these forms can be computed explicitely.
Then, it can be checked that not all of them are zero, and they do not belong to the 
families known by the work of Gekeler and Goss \cite{Ge, Go2}. For example, one computes easily:
$$D_{q^2 - q}h = \frac{g^{q - 1}h^q}{[1]^{q - 1}},\quad D_q(\Delta)=\frac{gh^q}{[1]},$$
so that $g^{q-1}h^q$ and $gh^q$ are Hecke eigenforms. 
\subsection{An algorithm to compute higher derivatives.}\label{algorithm}

Let $\delta=(\delta_n)_{n\geq 0}$ be a higher derivation on a $C$-algebra $F$.
Then, $X$ being an indeterminate over $F$, the map (Taylor's homomorphism)
$${\mathcal T}_{X}^\delta:F\rightarrow F[[X]]$$
defined by $${\mathcal T}_{X}^\delta(x)=\sum_{n=0}^\infty(\delta_nx)X^n$$ is a
$C$-algebra homomorphism \cite[Section 27]{MH}.

Over $F[[X]]$ there also is the iterative higher derivation $\delta'=(\delta'_n)_{n\geq 0}$
uniquely determined by
\begin{equation}\label{eqdprim}
\delta'_n(fX^i)=\binom{i}{n}fX^{i-n}
\end{equation}
for $f\in F$ and $n,i\geq 0$.
One checks that $(\delta_n)_{n\geq 0}$ is iterative if and only if, over $F$,
\begin{equation}\label{equivalent_to_iterative}
{\mathcal T}_X^\delta\circ \delta_n=\delta'_n\circ {\mathcal T}_X^\delta,\quad n\geq 0.\end{equation}
Indeed, this condition is equivalent to the commutativity of the diagram on p. 209 of \cite{MH}.


From now on we work with $F=C(E,g,h)$ and $\delta=D=(D_n)_{n\geq 0}$.
We also write, to ease notations, $\mathcal{T}_{X}=\mathcal{T}_{X}^D$.
We will also look at the fraction field $F(X)$ of $F[X]$ as embedded in $F((X))$.
For example, the expression $1/(1-EX)$ represents the formal series $1+EX+E^2X^2+\cdots\in F[[X]]$.

Let $f$ be an element of $F$.
In the following, we will make use of the polynomials
$${\mathcal T}_{X,k}(f):=\sum_{i=0}^{q^k-1}(D_if)X^i\in F[X],$$ so that, in $F[[X]]$, we have the following congruence modulo the ideal $(X^{q^k})$:
$${\mathcal T}_X(f)\equiv {\mathcal T}_{X,k}(f)\pmod{(X^{q^k})}.$$
These polynomials provide approximations to arbitrary order for the formal series $\mathcal{T}_X(f)$ as
$\mathcal{T}_X(f)=\lim_{n\rightarrow\infty}{\mathcal T}_{X,n}(f)$ (limit for the $X$-adic metric).

The map ${\mathcal T}_{X,k}:F\rightarrow F[X]$ is not a $C$-algebra homomorphism itself, but
induces a $C$-algebra homomorphism:
$${\mathcal T}_{X,k}:F\rightarrow \frac{F[[X]]}{(X^{q^k})}=\frac{F[X]}{(X^{q^k})}.$$

The following identities and congruences will prove to be useful
(the product being equal to $1$ when the indexing set is empty):
\begin{eqnarray}
\mathcal{T}_{X,r+s}(f^{q^s})&=&\mathcal{T}_{X,r}(f)^{q^s}\label{easy1}\\
\mathcal{T}_{X,r}(f^{-1})&\equiv& f^{-q^{r}}\prod_{i=0}^{r-1}\mathcal{T}_{X,r-i}(f^{q-1})^{q^i}\pmod{(X^{q^{r}})}\label{easy2}.
\end{eqnarray}
Equality (\ref{easy1}) holds for $f\in F$ and $r,s\geq 0$ and its validity is easy to check.
Congruence (\ref{easy2}) holds for $f\in F^\times$ and $r\geq 0$ and can be proved as follows.
Since $\mathcal{T}_{X}(f)^{q^{r}}\equiv f^{q^{r}}\pmod{(X^{q^{r}})}$, we have
$$\mathcal{T}_{X,r}(f^{q-1})^{1+q+\cdots+q^{r-1}}\mathcal{T}_{X,r}(f)\equiv f^{q^r}\pmod{(X^{q^r})},$$ yielding the desired congruence.

\begin{Proposition}\label{proposition_algorithm} Let $r,s\geq 0$ be integers. The following 
congruence holds, modulo the ideal $(X^{q^{r+s+1}})$ of $F[[X]]$:
\begin{eqnarray}\label{formulaproposition}
\lefteqn{{\mathcal T}_{X,r+s+1}(g_s)\equiv}\label{r+s+1}\\
&\equiv &[s+1]^{-1}{\mathcal T}_{X,r+1}(\Delta^{-1})^{q^s}
({\mathcal T}_{X,r}(g)^{q^{s+1}}{\mathcal T}_{X,r+s+1}(g_{s+1})
-{\mathcal T}_{X,r+s+1}(g_{s+2}))\nonumber.\end{eqnarray}
\end{Proposition}

\noindent\emph{Proof.} We appeal to the formula (\ref{eq:g0g1gk}), which is equivalent to
$$g_s=[s+1]^{-1}\Delta^{-q^s}(g^{q^{s+1}}g_{s+1}-g_{s+2}),\quad s\geq 0.$$
Congruence (\ref{r+s+1}) is obtained applying the $C$-algebra homomorphism ${\mathcal T}_X$  to both sides
of the latter identity (formulas  (\ref{easy1}) and  (\ref{easy2}) can help):
\begin{equation}\label{psiXgs}
{\mathcal T}_X(g_s)=[s+1]^{-1}{\mathcal T}_X(\Delta)^{-q^s}({\mathcal T}_X(g)^{q^{s+1}}{\mathcal T}_X(g_{s+1})-{\mathcal T}_X(g_{s+2})),\quad s\geq 0,
\end{equation}
and then reducing modulo the ideal $(X^{q^{r+s+1}})$ of $F[[X]]$.\CVD

We now give an algorithm for the explicit computation of
$\mathcal{T}_X(g_1),\mathcal{T}_X(g_2),\ldots$ and
$\mathcal{T}_X(\Delta)$, for which Proposition \ref{proposition_algorithm}
and Lemma \ref{lemme:2ek} are the key tools.

The best way to describe our algorithm, naturally presented as an
induction process, is to begin by giving the detail
of the explicit computation of its first steps.
We start with the explicit computation of the polynomials
$\mathcal{T}_{X,1}(\Delta)$ and $\mathcal{T}_{X,s+1}(g_s)$ ($s\geq 1$)
(first approximation). Then, we proceed to the computation of
$\mathcal{T}_{X,2}(\Delta)$ and of a representative of
$\mathcal{T}_X(g_s)$ modulo the ideal $(X^{q^{s+1}+2})$
(second approximation). Unfortunately, the full computation of the
polynomials $\mathcal{T}_{X,s+2}(g_s)$ would be too long to present
in this article.
 
These explicit computations will prepare the reader for the general
process which generates explicit expressions for the polynomials
$\mathcal{T}_{X,r+s+1}(g_s)$ ($s\geq 1$) and $\mathcal{T}_{X,r+1}(\Delta)$
for all $r\geq 1$; he or she will then be ready to understand the
algorithm. At the same time, the accomplished explicit computations
are used for different purposes in several parts of this paper.

In \cite[Theorem 4.1]{BP} we have computed $D_1f,D_qf,D_{q^2}f\in C[E,g,h]$
with $f\in\{E,g,h\}$, applying the classical technique consisting
in solving linear equations in $C$-vector spaces of modular forms with
prescribed order of vanishing at infinity. This method of computation
can be pushed beyond to compute also $D_{q^3}f,D_{q^4}f,\ldots$ but then
it requires that one first computes the coefficients
of the $u$-expansions of $E,g,h$ with Gekeler's algorithm of \cite{Ge},
before entering the linear algebra part.
However, the computation of the $u$-expansions of $E,g,h$ is not an
easy matter, since it also needs computation of the so-called Goss
polynomials, a task that usually generates large computations.

The algorithm we give here is of a different nature and can be considered
as a variant of Gekeler's techniques of computing Goss polynomials (cf. Proof of Lemma \ref{inseparable}), and using the
recurrence relations (\ref{eq:g0g1gk}) to yield $u$-expansions of
modular forms. Our algorithm is easier to use, compared to the methods
introduced in \cite{BP} because it does not need any preliminary
computation of $u$-expansions.

\subsubsection{First approximation.}

We know that:
\begin{eqnarray}
\mathcal{T}_{X,s}(g_s)&=&g_s+x_sX,\quad s\geq 1,\label{eqTxsgs}\\
\mathcal{T}_{X,1}(\Delta)&=&\Delta(1-EX)\label{TX1Delta}.
\end{eqnarray} 
The first formula follows easily from Lemma \ref{lemme:2ek}, while
the second can be obtained  by using \cite[Theorem 4.1, (iii)]{BP}
or applying Lemma \ref{lemme:2ek} to $g_1,g_2$
and then use the formula 
\begin{equation}
\Delta=[1]^{-1}(g_1^{q+1}-g_2),\label{eqDelta}
\end{equation}
easily deduced from (\ref{eq:g0g1gk}).

From (\ref{TX1Delta}) we obtain:
\begin{eqnarray*}
{\mathcal T}_{X,1}(\Delta^{-1})&\equiv&\Delta^{-1}\frac{1}{1-EX}\pmod{(X^q)}\\
&=&\Delta^{-1}(1+EX+E^2X^2+\cdots+E^{q-1}X^{q-1}).\end{eqnarray*}
Substituting (\ref{eqTxsgs}) in (\ref{formulaproposition}) 
we obtain the following congruence modulo the ideal $(X^{q^{s+1}})$
for ${\mathcal T}_{X,s+1}(g_s)$, where we suppose that $s\geq 1$:
\begin{eqnarray}
\lefteqn{{\mathcal T}_{X,s+1}(g_s)\equiv}\nonumber\\
&\equiv&[s+1]^{-1}\Delta^{-q^s}\frac{1}{1-E^{q^s}X^{q^s}}((g^{q^{s+1}}
+x_1^{q^{s+1}}X^{q^{s+1}})(g_{s+1}+x_{s+1}X)-\nonumber\\
& &-(g_{s+2}+x_{s+2}X))\nonumber\\
&\equiv&\frac{g_s+x_sX}{1-E^{q^s}X^{q^s}}\label{Ts+1gs}.
\end{eqnarray}
We thus obtain the following formula for ${\mathcal T}_{X,s+1}(g_s)$:
$${\mathcal T}_{X,s+1}(g_s)=g_s+x_sX+E^{q^s}X^{q^s}(g_s+x_sX)+\cdots+E^{(q-1)q^s}X^{(q-1)q^s}(g_s+x_sX).$$
This in turn allows, inserting the congruence (\ref{Ts+1gs}) for $s=1,2$ in (\ref{eqDelta}), to compute ${\mathcal T}_{X,2}(\Delta)$. Here is the formula we find (congruences modulo the ideal $(X^{q^2})$):
\begin{eqnarray}
\lefteqn{{\mathcal T}_{X,2}(\Delta)\equiv [1]^{-1}(\mathcal{T}_{X,2}(g_1)^q\mathcal{T}_{X,2}(g_1)-\mathcal{T}_{X,2}(g_2))}\nonumber\\ &\equiv&[1]^{-1}\left(\frac{(g+x_1X)^{q+1}}{(1-E^{q}X^{q})^{q+1}}-(g_2+x_2X)\right)\nonumber\\
&\equiv&\Delta\left(1-EX+\left(\frac{gh}{[1]}-E^q\right)X^q+\right.\nonumber\\
& &\left.\left(E^{q+1}-\frac{Egh}{[1]}-\frac{h^2}{[1]}\right)X^{q+1}\right)(1-E^qX^q)^{-1}.\label{TX2Delta1}
\end{eqnarray}
To check the last congruence, the reader can make use of the tables of Section \ref{tables} and identity (\ref{eqDelta}).

\subsubsection{Second approximation.}\label{second_approx}
By using (\ref{r+s+1}) and (\ref{Ts+1gs}), we get, for all $s\geq 1$,
the following congruences modulo the ideal $(X^{q^{s+2}})$:
\begin{eqnarray*}
\lefteqn{[s+1]{\mathcal T}_{X,s+2}(g_s)\equiv}\\
&\equiv&{\mathcal T}_{X,2}(\Delta^{-1})^{q^s}({\mathcal T}_{X,1}(g)^{q^{s+1}}{\mathcal T}_{X,s+2}(g_{s+1})-{\mathcal T}_{X,s+2}(g_{s+2}))\\
&\equiv&{\mathcal T}_{X,2}(\Delta^{-1})^{q^s}\left((g^{q^{s+1}}+x_1^{q^{s+1}}X^{q^{s+1}})\frac{g_{s+1}+x_{s+1}X}{1-E^{q^{s+1}}X^{q^{s+1}}}-(g_{s+2}+x_{s+2}X)\right).
 \end{eqnarray*}
 
The explicit computation of ${\mathcal T}_{X,2}(\Delta^{-1})$ is possible with Formula (\ref{easy2}).
 Unfortunately, it is rather complicated to handle, so we 
 limit ourselves to its determination modulo the ideal $(X^{q+1})$. This is why we do not fully compute $\mathcal{T}_{X,s+2}(g_s)$ in this text.
 
 Since ${\mathcal T}_{X,2}(\Delta)\equiv\Delta(1-EX+([1]^{-1}gh-E^q)X^q)\pmod{(X^{q+1})}$ we compute easily:
 \begin{equation}
 {\mathcal T}_{X,2}(\Delta^{-1})\equiv\Delta^{-1}\left(\sum_{i=0}^{q-1}E^iX^i+\left(E^q-\frac{gh}{[1]}\right)X^q\right)\pmod{(X^{q+1})}.\label{TX2Delta}
 \end{equation}
We reduce the $q^s$-th power of the polynomial (\ref{TX2Delta}) modulo the ideal $(X^{q^{s+1}+2})$.
Looking at the formula (\ref{eq:g0g1gk}) in the form $$\Delta^{q^{s}}g_s=[s+1]^{-1}(g^{q^{s+1}}g_{s+1}-g_{s+2})$$ and $x_1=-Eg-h$, we find the congruences (modulo $(X^{q^{s+1}+2})$):
\begin{eqnarray*}
\lefteqn{[s+1]{\mathcal T}_{X}(g_s)\equiv}\\
&\equiv &{\mathcal T}_{X,2}(\Delta^{-1})^{q^s}((g^{q^{s+1}}+x_1^{q^{s+1}}X^{q^{s+1}})(g_{s+1}+x_{s+1}X)(1+E^{q^{s+1}}X^{q^{s+1}})-\\ & &-(g_{s+2}+x_{s+2}X))\\
&\equiv&{\mathcal T}_{X,2}(\Delta^{-1})^{q^s}(g^{q^{s+1}}g_{s+1}-g_{s+2}+(g^{q^{s+1}}x_{s+1}-x_{s+2})X+\\ & &
+(x_1^{q^{s+1}}g_{s+1}+E^{q^{s+1}}g_{s+1}g^{q^{s+1}})X^{q^{s+1}}-h^{q^{s+1}}x_{s+1}X^{q^{s+1}+1})\\
&\equiv&{\mathcal T}_{X,2}(\Delta^{-1})^{q^s}(\Delta^{q^s}[s+1]g_s+\Delta^{q^s}[s+1]x_sX-\\
& & -h^{q^{s+1}}g_{s+1}X^{q^{s+1}}-h^{q^{s+1}}x_{s+1}X^{q^{s+1}+1})\\ & &\\
&\equiv&\left(\sum_{i=0}^{q-1}E^{iq^s}X^{iq^s}+\left(E^{q^{s+1}}-\frac{g^{q^s}h^{q^s}}{[1]^{q^s}}\right)X^{q^{s+1}}\right)\times\\  & &
([s+1]g_s+[s+1]x_sX+h^{q^{s}}g_{s+1}X^{q^{s+1}}+h^{q^{s}}x_{s+1}X^{q^{s+1}+1})\\
& &\pmod{(X^{q^{s+1}+2})}.
\end{eqnarray*}
Therefore, we obtain:
\begin{Proposition}\label{propformulasgsxs} For $s\geq 1$, we have the congruences:
\begin{eqnarray}\label{xigs}
\lefteqn{{\mathcal T}_X(g_s)\equiv}\\ & \equiv&\frac{g_s+x_sX}{1-E^{q^s}X^{q^s}}-h^{q^s}\left(-\frac{g_{s+1}}{[s+1]}+\frac{g^{q^s}g_s}{[1]^{q^s}}\right)X^{q^{s+1}}-\nonumber\\ & &
-h^{q^s}\left(-\frac{x_{s+1}}{[s+1]}+\frac{g^{q^s}x_s}{[1]^{q^s}}\right)X^{q^{s+1}+1}\pmod{(X^{q^{s+1}+2})},\nonumber\\
\lefteqn{{\mathcal T}_X(x_s)\equiv}\label{xixs}\\ & \equiv&\frac{x_s}{1-E^{q^s}X^{q^s}}-h^{q^s}\left(-\frac{x_{s+1}}{[s+1]}+\frac{g^{q^s}x_s}{[1]^{q^s}}\right)X^{q^{s+1}}\pmod{(X^{q^{s+1}+1})}.\nonumber
\end{eqnarray}
In particular, $\epsilon_D(x_k)=(k+1)e$ for all $k\geq 0$.
\end{Proposition}
The second congruence of the proposition follows from the first and (\ref{equivalent_to_iterative}) because we have, for $s\geq 1$:
$$\mathcal{T}_{X}(x_s)=\mathcal{T}_X(D_1g_s)=\frac{\partial}{\partial X}\mathcal{T}_X(g_s).$$
The value of $\epsilon_D(x_k)$ is now easy to determine thanks to (\ref{xixs}). Indeed, the first term of the right-hand side 
of this identity tell that $\epsilon_D(x_k)\geq (k+1)e$. Equality follows by checking that the coefficient of $X^{q^{k+1}}$
is a polynomial of $C[E,g,h]$ which is coprime with $x_k$ by Lemma
\ref{xkirreducible}.\CVD

\begin{Remark}{\em 
Notice that, substituting $s=0$, the second formula agrees with \cite[(iv) and (vii) of Theorem 4.1]{BP}.
There is no simple explanation of this fact.}
\end{Remark}

\subsubsection{End of description of the algorithm\label{description_algorithm}.}

We shall now describe, in its generality, the algorithm, which uses induction on $r\geq 1$. Taking into account Proposition \ref{propformulasgsxs}
and formula (\ref{TX2Delta}), we can assume that for an integer $r\geq 1$ we have already computed $\mathcal{T}_{X,r+s}(g_s)$ (for all $s\geq 1$)
and $\mathcal{T}_{X,r+1}(\Delta)$, explicitly as polynomials of $C(E,g,h)[X]$.

Therefore, in (\ref{r+s+1}), the polynomials $${\mathcal T}_{X,r+s+1}(g_{s+2}),\quad
{\mathcal T}_{X,r}(g),\quad{\mathcal T}_{X,r+s+1}(g_{s+1})$$ are all known.
By using (\ref{easy2}), we compute $\mathcal{T}_{X,r+1}(\Delta^{-1})$, then we use (\ref{easy1}) to raise the expression to the $q^s$-th power;
we obtain the polynomials $\mathcal{T}_{X,r+s+1}(g_s)$, $s\geq 1$.

In particular, for $s=1,2$ we have an explicit expression of $\mathcal{T}_{X,r+2}(g_1)$
and $\mathcal{T}_{X,r+3}(g_2)$ and the degree $< q^{r+2}$ representative of the class of reduction modulo $(X^{q^{r+2}})$ of the polynomial
$$[1]^{-1}(\mathcal{T}_{X,r+2}(g_1)^{q+1}-\mathcal{T}_{X,r+3}(g_2))$$
is the polynomial ${\mathcal T}_{X,r+2}(\Delta)$ by (\ref{eqDelta}).

This ends the description of the algorithm.


\subsubsection{An algorithm which computes $D_nf$ for $f=E,g,h$, $n\geq 1$.}

It remains to explain how to compute the polynomials $\mathcal{T}_{X,r}(f)$ with $f\in\{E,g,h\}$ for $r\geq 0$.
First of all, $\mathcal{T}_{X,r}(g)$ and $\mathcal{T}_{X,r}(\Delta)$
are computed by the algorithm of section \ref{description_algorithm} (case $s=1$).

We have 
$\Delta=-h^{q-1}$. Hence, $\mathcal{T}_X(\Delta)=-{\mathcal T}_{X}(h)^{q-1}$.
By \cite[Proposition 3.6]{BP}, there exist two formal series
$$f_h=\sum_{m=0}^\infty c_mX^m,\quad f_\Delta=\sum_{m=0}^\infty d_mX^m\in F[[X]]$$
with $c_0=d_0=1$, such that
${\mathcal T}_X(h)=hf_h$ and $ {\mathcal T}_X(\Delta)=\Delta f_\Delta$.
We have $$f_h^{q-1}=f_\Delta,$$ so that the coefficients $c_n$ are uniquely determined by 
the following relations:
$$c_n=\ell_n-\sum_{i+j=n, i\not=n}c_id_j,\quad \text{ where }\ell_n=\left\{\begin{array}{ll}0&\text{ if }q\nmid n\\
c_{n/q}^q&\text{ if }q|n\end{array}\right.,\quad n\geq 0.$$
The coefficients $c_n$ can be computed by induction on $n\geq 0$ and this allows
to compute  the polynomials ${\mathcal T}_{X,r}(h)$.

Now, $D_1h=Eh$, so that ${\mathcal T}_X(D_1h)={\mathcal T}_X(E){\mathcal T}_X(h)$. But ${\mathcal T}_X(D_1h)=\delta'_1{\mathcal T}_X(h)$
because $D$ is iterative and we have (\ref{equivalent_to_iterative}), so that $${\mathcal T}_X(E)=(\delta'_1{\mathcal T}_X(h)){\mathcal T}_X(h)^{-1}$$
where $(\delta'_n)_{n\in\NN}$ is the iterative derivation defined in (\ref{eqdprim}). Therefore,
\begin{equation}\mathcal{T}_{X,r}(E)\equiv\left(\frac{\partial}{\partial X}{\mathcal T}_{X,r}(h)\right){\mathcal T}_{X,r}(h^{-1})\pmod{(X^{q^r})},\quad r\geq 0.
\label{Efromh}\end{equation}
The computation of the polynomials ${\mathcal T}_{X,r}(h^{-1})$ can be made 
with (\ref{easy2}).

In fact, since $D_1(\Delta)=-E\Delta$, the computation of the sequence ${\mathcal T}_{X,r}(E)$
can be achieved avoiding the use of $h$, by using (\ref{equivalent_to_iterative}).\CVD

\begin{Remark}{\em  By Theorem 1 and Proposition 3.1 of \cite{BP},
we also know that all the fractions $D_nE,D_ng,D_nh\in F$ belong to $\widetilde{M}=C[E,g,h]$.
In fact, we have $\mathcal{T}_X(f)\in \widetilde{M}[[X]]$ for all $f\in\widetilde{M}$.
This property seems not to follow from Proposition \ref{proposition_algorithm}.
}\end{Remark}

\medskip

\noindent\emph{Acknowledgement.} The Referee rapidly provided us
with a detailed report which helped us to improve the present text. The Referee highlighted a mistake in a previous proof we gave of Theorem \ref{theorem:xkbis}. It is in the subsequent correction we have 
proposed, that we have appealed to \cite{GG}.
For all this we would like to express our gratitude.

\vspace{15pt}

\noindent{\small Vincent Bosser,\\
\indent L.M.N.O., Universit\'e de Caen,\\
\indent Campus II - Boulevard Mar\'echal Juin,\\
\indent BP 5186 - F14032 Caen Cedex, France.\\
\indent E-mail: {\tt bosser@math.unicaen.fr}}

\vspace{15pt}

\noindent{\small Federico Pellarin,\\
\indent Laboratoire La.M.U.S.E., Universit\'e de Saint-Etienne,\\
\indent Facult\'e de Sciences - 23, rue du Dr. P. Michelon,\\
\indent 42023 Saint-Etienne Cedex, France.\\
\indent E-mail: {\tt federico.pellarin@univ-st-etienne.fr}}

\end{document}